\documentclass[12pt]{amsart} 
\usepackage{epsf}

\def\ie{{\it i.e.}}

\def\cf{{\it cf.}}
\def\viceversa{{\it vice versa}}
\def\fg{fundamental group}
\def\fgs{fundamental groups}

\def\arr{{arrangement}}
\def\arrs{{arrangements}}
\def\Empty{{\small{no\ multiple\ points}}}

\newcommand\sg[1]{{\left<{#1}\right>}}

\newtheorem{thm}{Theorem}[section]
\newtheorem{prop}[thm]{Proposition}
\newtheorem{cor}[thm]{Corollary}
\newtheorem{lem}[thm]{Lemma}

\newtheorem{defn}[thm]{Definition}

\newtheorem{rem}[thm]{Remark}

\newcommand{\barr}{\begin{array}}
\newcommand{\earr}{\end{array}}
\newcommand{\btab}{\begin{tabular}}
\newcommand{\etab}{\end{tabular}}
\newcommand{\beq}{\begin{equation}}
\newcommand{\eeq}{\end{equation}}
\newcommand{\bea}{\begin{eqnarray*}}
\newcommand{\eea}{\end{eqnarray*}}
\newcommand{\bce}{\begin{center}}
\newcommand{\ece}{\end{center}}
\newcommand{\bpi}{\begin{picture}}
\newcommand{\epi}{\end{picture}}
\long\def\forget#1\forgotten{}
\newcommand{\bsl}{\begin{slide}{}}
\newcommand{\esl}{\end{slide}}

\newcommand{\bib}{thebibliography}

\newcommand{\sbs}{\subset}

\newcommand\set[1]{{\{{#1}\}}}

\def\suchthat{{\,|\,}}
\newcommand{\halfequiv}{\equiv_{\frac{1}{2}}}

\newcommand\software[1]{{\textsf{#1}}}
\newcommand\Lpair[2]{{\left[{#1},{#2}\right]}}

\newcommand{\al}{\alpha}
\newcommand{\be}{\beta}
\newcommand{\ga}{\gamma}

\newcommand{\si}{\sigma}

\newcommand{\Ga}{\Gamma}

\newcommand{\ba}{{{\bf {a}}}}
\newcommand{\bb}{{\bf b}}

\newcommand{\C}{{{\mathbb{C}}}}

\newcommand{\R}{{\mathbb R}}
\newcommand{\Z}{{\mathbb Z}}
\newcommand{\F}{{\mathbb F}}
\newcommand{\PP}{{\mathbb P}}
\newcommand{\cal}{\mathcal}

\newcommand{\cL}{{\cal L}}

\newcommand{\LP}{Lefschetz pair}
\newcommand{\LPs}{Lefschetz pairs}

\def\Wmodequiv{{W_S/\!\equiv}}
\newcommand\tru[3][c+]{{{#1}\triangle_{#2}^{(#3)}}}
\newcommand\trd[3][c+]{{{#1}\nabla_{#2}^{(#3)}}}

\newcommand\figs[1]{#1}
\newcommand\arel{{\stackrel{\Delta}{=}}}
\newcommand\defin[1]{{\it{#1}}}
\newcommand\FIGURE[4][]{{\begin{figure}[!h]\epsfysize=#3 {\epsfbox{\figs{#2}}}\caption{#1}\label{#4}\end{figure}}}
\newcommand\FIGUREx[4][]{{\begin{figure}[!h]\epsfxsize=#3 {\epsfbox{\figs{#2}}}\caption{#1}\label{#4}\end{figure}}}
\newcommand\tabfig[2]{{$\stackrel{\ }{\epsfysize=#1 \epsfbox{\figs{#2}}}$}}

\newif\iffurther
\furtherfalse

\begin{document}

\title[Classification of Line Arrangements]%
    {$\pi_1$-classification of real arrangements with up to eight lines}

\author[David Garber, Mina Teicher]{David Garber$^1$, Mina Teicher$^1$}
\email{\set{garber,teicher}@macs.biu.ac.il}
\address{
Dept. of Math. and CS,
Bar-Ilan University,
Ramat-Gan 52900, Israel
}
\author[Uzi Vishne]{Uzi Vishne$^2$}
\email{vishne@math.huji.ac.il}
\address{
Einstein Institute of Mathematics, the Hebrew University,
Jerusalem 91904, Israel
}

\forget
\author
\author
{David Garber$^1$ \and Mina Teicher$^1$ \and Uzi Vishne$^2$}

\address{
Dept. of Math and CS. \\
Bar-Ilan University \\
Ramat-Gan 52900, Israel\\
\\
Einstein Inst. of Math.\\
Hebrew University \\
Jerusalem 91904, Israel
}
\forgotten

\date{\today}

\stepcounter{footnote}
\footnotetext{Partially supported by The Israel Science Foundation
(Center of Excellence Program), and by the Emmy Noether Institute for
Mathematics and by the Minerva Foundation (Germany).
The research was done during the Ph.D. studies of David Garber,
under the supervision of Prof. Mina Teicher.}

\stepcounter{footnote}
\footnotetext{Partially supported by the Edmund Landau Center for
Research in Mathematical Analysis and Related Subjects.}

\begin{abstract}
One of the open questions in the geometry of line arrangements
is to what extent
does the incidence lattice of an arrangement determine its fundamental group.
Line arrangements of up to 6 lines were recently classified by K.M. Fan
\cite{Fa2}, and it turns out that the incidence lattice of
such arrangements determines the projective fundamental group.
We use actions on the set of wiring diagrams, introduced in \cite{GTV},
to classify real arrangements of up to 8 lines.
In particular, we show that the incidence lattice of such arrangements
determines both the affine and the projective fundamental groups.

\end{abstract}

\maketitle



\section{Introduction}
A \defin{line arrangement } in $\C^2$ is a union of finitely many copies of $\C^1$.
An arrangement is equipped with several invariants, probably
the most important are the fundamental groups of the complement in $\C^2$
and in $\C\PP^2$. These are  called the \defin{affine} and
the \defin{projective \fg} of the arrangement, respectively.

A more combinatorial invariant is the incidence lattice of the
arrangement. It is not known whether or not the lattice
determines the fundamental groups \cite{CS}.

In 1997, Fan \cite{Fa2} studied the projective fundamental group of arrangements with
small number of lines. Using classification of the arrangements with $6$ lines,
he was able to show that if there are up to $6$ lines, the incidence lattice determines
the projective fundamental group.

An arrangement is called \defin{real } if the defining equations of its lines 
can be written with real coefficients, and \defin{complex} otherwise.

Rybnikov, in an unpublished work \cite{Ry},
presents two complex arrangements with $13$ lines based on 
MacLane configuration,
which have the same lattice, but different projective fundamental groups.
This example was not thoroughly checked, and in any case the question for real arrangements remains open.

To a real arrangement, one can associate combinatorial objects: the wiring
diagram and its associated list
of Lefschetz pairs. They are not invariants, as they depend
on the choice of a guiding generic line. Still, these objects turn out
to be useful tools in the study of fundamental groups of line arrangements,
\cf\ \cite{MoTe1}, \cite{GaTe}, \cite{CS}.

\medskip

The purpose of this paper is twofold. We extend Fan's result to
arrangements of up to 8 lines (for the affine and projective
fundamental groups), and at the same time
present algorithms developed for this classification, which can be also used in other
classification problems.
While it is difficult to enumerate geometric objects as lines
(the defining coefficients
have infinitely many possible values), it is possible, in principle,
to enumerate the wiring diagrams (or the lists of Lefschetz pairs) induced by them.
However, many wiring diagrams are induced by the same line arrangement
(because they depend on the guiding line), and we do not want to waste time studying
various diagrams induced by the same arrangement.

To this end, we use equivalence relations and actions on
the set of wiring diagrams which were introduced in \cite{GTV}.
We proved there that the equivalence relations and the actions preserve our two
main invariants: the incidence lattice and the affine and projective
fundamental groups.
These results enable us to list all the incidence lattices of arrangements
up to 8 lines.
They also reduce the number of comparisons of pairs of fundamental
groups needed to test the conjecture from over $20$ millions to
about $200$ pairs. 
Eventually, we conclude that for real
arrangements with up to $8$ lines the incidence lattice
determines both the affine and the projective fundamental groups.

\medskip

The paper is organized as follows. In Section \ref{defs} we briefly recall
some combinatorial objects related  to a line arrangement:
the wiring diagram and the associated list of Lefschetz pairs, the incidence lattice
and the signature.

In section \ref{sec3} we discuss the affine and projective fundamental groups of some special
types of line arrangements. This is later used to avoid some cases in the classification.

In section \ref{actions}, we briefly survey equivalence relations
and actions on the set of wiring diagrams of a given signature.
These were introduced and studied in \cite{GTV}.

In Sections \ref{6lines}--\ref{8lines} we tabulate all the real
line arrangements of up to 8 lines according to their incidence
lattices (some details on the arrangements of 6 lines were
already given by Fan \cite{Fa2}). The classification is achieved
by classifying the appropriate lists of Lefschetz pairs, and
requires considerable computer work. Details on the algorithms
used are given in Section \ref{alg}. Our main interest in this
paper is the algorithms and the resulting classification. For
more details (and proofs) on the actions and relations defined
here, the reader is directed to \cite{GTV}.

\section{Combinatorial objects related to line arrangements} \label{defs}

We briefly recall some of the combinatorial objects and construction related to real
line arrangements. For more details, see \cite{GP} and \cite{GTV}.

\subsection{Wiring diagrams and Lefschetz pairs}

Given a real line arrangement in $\C^2$, its intersection with the natural copy
of $\R^2$ in $\C^2$ is an arrangement of lines in the real plane.

To an arrangement of $\ell$ lines in $\R^2$ one can associate a
\defin{wiring diagram} \cite{Go}, which holds the information
on the intersection points and the relative position of their
projection on a generic line. A wiring diagram is a collection of
$\ell$ pseudolines (where a \defin{pseudoline} in $\R^2$ is a
union of segments and rays, homeomorphic to $\R$). A line $L$ is
\defin{generic} with respect to a line arrangement $\cL$, if it
avoids all the intersection points of the arrangement, and the
projections of intersection points on $L$ do not overlap. Fixing
a generic line $L$, the wiring diagram induced by $\cL$ (with
respect to $L$) contains $\ell$ pseudolines, and is constructed
as follows: start at the '$-\infty$' side of $L$ with $\ell$
parallel rays, and for every projection of an intersection point,
make the corresponding switch in the rays.

To a wiring diagram, one can associate a list of \defin{Lefschetz
pairs}. Any pair of this list corresponds to one of the
intersection points, and holds the first and the last indices of
the pseudolines intersected at that point, numerated locally near
the intersection point, from bottom to top. An example is given in
Figure \ref{fig1}.

\FIGURE{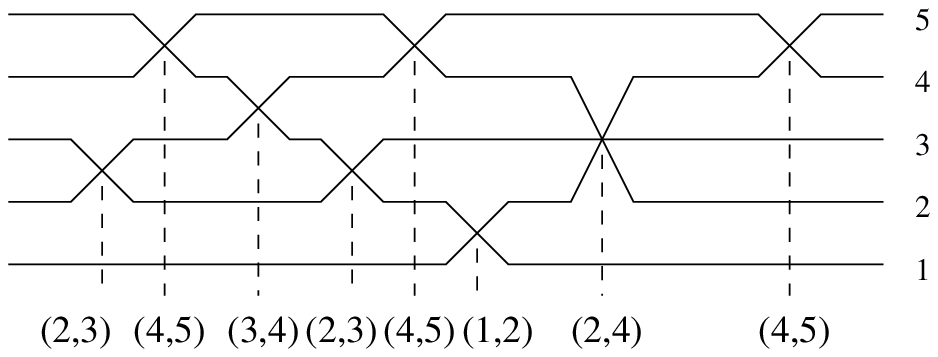}{4cm}{fig1}

To every list of Lefschetz pairs there corresponds a wiring
diagram, which is constructed by the reverse procedure. It should be noted that in
the wiring diagram  only the order
of the intersection points is relevant, and not their distances.

\subsection{The incidence lattice}
Let $\cL = \{ L_1,L_2, \cdots, L_{\ell} \}$ be a line arrangement. By ${\rm Lat}(\cL)$
we denote the partially-ordered set of non-empty intersections of the $L_i$, ordered by
inclusion (see \cite{OT}). We include the whole plane and the empty set
in ${\rm Lat}(\cL)$, so that it becomes a lattice of height $3$:
the intersection points have height $1$, and the lines height $2$.

Equivalently, the lines and points form a bipartite graph, in which a line $L_i$
and a point $p_j$ are connected iff $p_j \in L_i$.


\subsection{Lattice isomorphism and the fundamental group}


As mentioned in the introduction, it was conjectured by Cohen and
Suciu \cite{CS} that for real \arrs\ the incidence lattice determines the
(affine and projective) \fgs.
There is an unpulished example of Rybnikov \cite{Ry}
which shows that this is not the
case for complex \arrs\ (at least for the projective \fg).

A rather weak version of that conjecture was recently proved by
Cordovil \cite{Co}. He introduces a new order (called
\defin{shellability order}) on the lines of an arrangement. Then,
he shows that if there is a bijection of two arrangements which
induces an isomorphism of the associated lattices, preserving the
shellability order, then their affine fundamental groups
are isomorphic. 

\subsection{The signature}\label{signature}
Another characteristic of a line arrangement is the so-called
\defin{signature} - the number of lines meeting in every
intersection point.

\begin{defn}
The {\rm signature} of a line arrangement $\cL$ is $[2^{n_2}3^{n_3}\dots]$,
where $n_k$ is the number of points in which $k$ lines intersect.
We make the agreement to  omit every $c^0$.
\end{defn}

For example, an arrangement which consists of five lines, four of which
intersect at one point, has signature $[2^4 4^1]$. The signature is obviously
determined by the incidence lattice.

The definition of the signature applies to any wiring diagram (and thus also to the
associated list of Lefschetz pairs).

A point in which only two lines  intersect is called \defin{simple}, so that
$n_2$ is the number of simple points in the \arr. The other points are
called \defin{multiple}.

Note that different incidence lattices may still have the same signature.
For example, both the
following two arrangements have the same signature $[2^9 3^2]$,
but their incidence lattices are not isomorphic
(the two multiple points are connected in the first arrangement but not in the
second):

\FIGURE{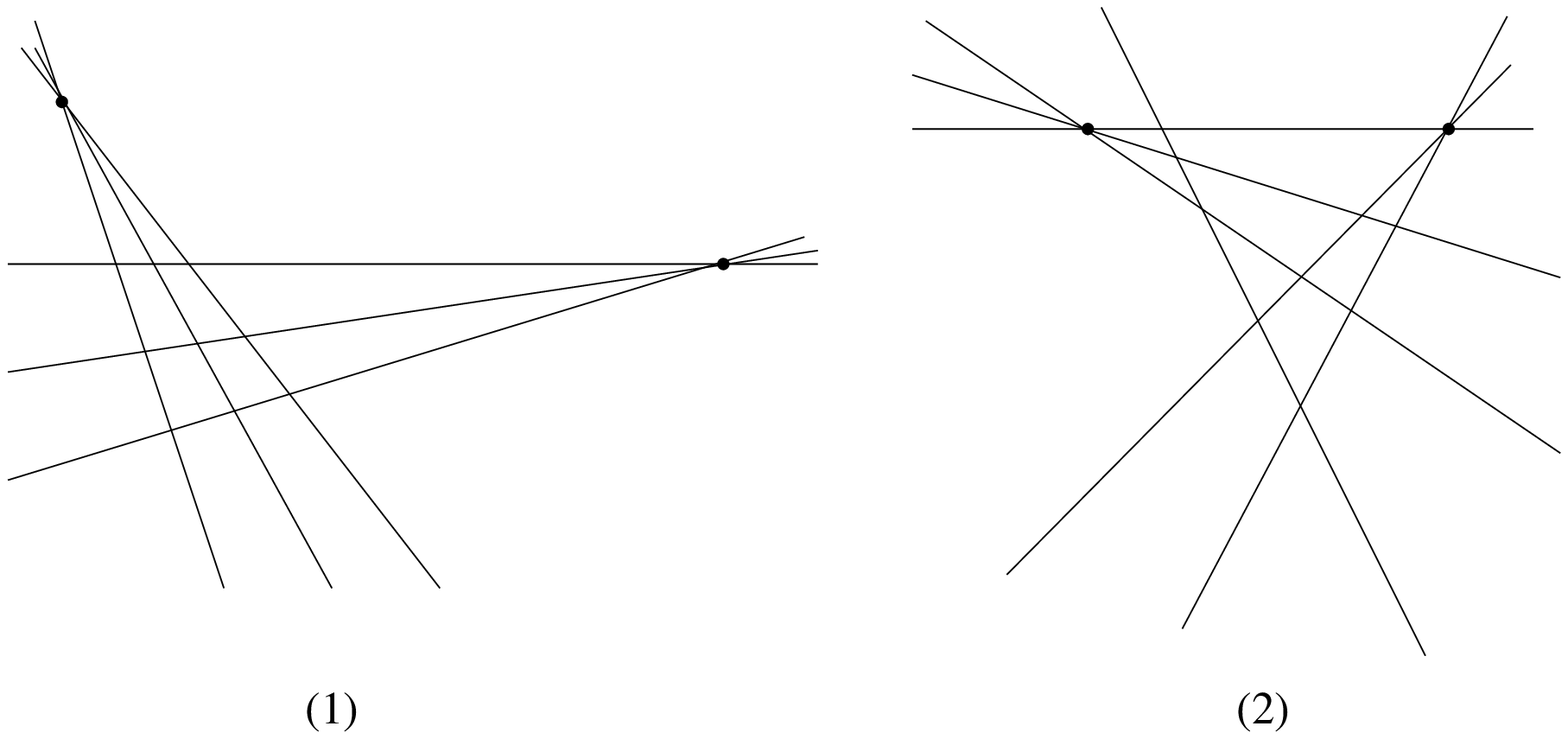}{4cm}{fig2}

\begin{rem}
The affine \fgs\ of both arrangements are isomorphic to
$\F_2\oplus \F_2\oplus \Z^2$, and the projective \fgs\ are isomorphic to
$\F_2\oplus \F_2 \oplus \Z$, \cite{GaTe}.
Here $\F _2$ is the free group on two generators.
\end{rem}

It is not difficult to check that the example above is minimal:
for less than $6$ lines, the signature of an arrangement does determine the
incidence lattice.

We say that a wiring diagram has the \defin{unique intersection property}, if every two wires
intersect exactly once. Obviously, a wiring diagram induced by a line arrangement has
this property.

\begin{lem} \label{UIP}
A wiring diagram has the unique intersection property if and only if the following
conditions hold:
\begin{enumerate}
\item Its signature $S=[2^{n_2} 3^{n_3} \cdots ]$ satisfies
\begin{equation}\label{SUIP}
\sum _{k\geq 2} n_k{{k} \choose {2}} = {{\ell} \choose {2}},
\end{equation}
where $\ell$ is the number of the wires in the diagram, and
\item The order of the wires at $+ \infty$ is the opposite of that at $- \infty$.
\end{enumerate}
\end{lem}

\begin{proof}
The unique intersection property implies that the total number of
intersection points is ${{\ell}\choose {2}}$. Moreover, if a wire
$x$ was below a wire $y$ at $+ \infty$ it will be above it from
their intersection point on.

Conversely, if every two pseudolines changed order, then we have at least
${{\ell}\choose {2}}$ intersection points, and since this is the total number, every two
lines intersect exactly once.
\end{proof}

Note that not every signature satisfying Equation
 (\ref{SUIP}) can arise, an easy counterexample
is given in the following remark:
\begin{rem}
There is no wiring diagram (with $\ell = 4$ pseudolines) whose signature
is $[3^2]$ (for if three lines intersect at one point,
the remaining line must either join them at that point, or be a simple line,
resulting in signature $[4^1]$ or $[2^3 3^1]$).
\end{rem}

\section{The fundamental group of a line arrangement}\label{sec3}

In this section we discuss the fundamental groups
of some special types  of real line arrangements.
This will be later used to analyse some special signatures.

\subsection{Arrangements with no cycles of multiple points}
The graph $G(\cL)$ of multiple points of a
line arrangement $\cL$ lies on the arrangement. It consists of the multiple
points of $\cL$, with the lines on which there are at least two multiple
points.

Fan has proved the following result:

\begin{thm}[{\cite[Theorem 3.2]{Fa2}}] \label{Fan}
If $G(\cL)$ has no cycles, then the projective fundamental group of $\cL$ is
$$\pi_1 (\C\PP^2 -\cL) \cong \F _{m_1 -1} \oplus \cdots \oplus \F _{m_k -1}
\oplus \Z ^{\ell-({\sum_{i=1}^k} (m_i -1))-1},$$
where $m_1, \cdots , m_k$ are the multiplicities of the multiple intersection points
in $\cL$ and $\ell$ is the number of lines.
\end{thm}

A somewhat weaker result can be proved, using braid monodromy
techniques and the van Kampen theorem, for the affine fundamental
group:
\begin{thm}[{\cite[Theorem 5.3]{GaTe}}]\label{Fan-}
If every connected component of $G(\cL)$ consists of vertices lying on
a straight line, then the affine fundamental group is:
$$\pi _1 (\C ^2 -\cL) \cong \F _{m_1 -1} \oplus \cdots \oplus \F _{m_k -1}
\oplus \Z ^{\ell-({\sum_{i=1}^k} (m_i -1))}$$
where $m_1, \cdots , m_k$ and $\ell$ are as in Theorem \ref{Fan}.
\end{thm}

It immediately follows that if  $\cL$ has at most two
multiple points, then the (affine and projective)
fundamental groups depend only on the signature of $\cL$.

For other results on the topology
of the complement of real arrangements, see \cite{Ra},
\cite{Sa}, \cite{Ar}, or \cite{OrSo}.

\subsection{Arrangements with simple lines}
We call a line $L \in \cL$ \defin{simple} if it has only simple
intersection points. Note that if $\cL=\cL'
\cup L$ where $L$ is a simple line in $\cL$, then the incidence
lattice of $\cL$ is determined by that of $\cL'$.

The following result of Oka and Sakamoto allows us to reduce (for affine
fundamental groups) to the case of no simple lines:

\begin{prop}[\cite{OkSa}] \label{OS}
Let $C_1$ and $C_2$ be
algebraic plane curves in $\C ^2$ of degrees $d_1$,$d_2$ respectively.
Assume that the intersection $C_1 \cap C_2$
consists of $d_1 d_2$ distinct points. Then:
$$\pi _1 (\C ^2 - (C_1 \cup C_2)) \cong \pi _1 (\C ^2 -C_1) \oplus \pi _1 (\C ^2 -C_2).$$
\end{prop}

In particular, if $\cL=\cL' \cup L$ where $L$ is a simple line in $\cL$, then:
$$\pi _1(\C ^2 -\cL) \cong \pi _1(\C ^2 -\cL')\oplus \pi_1(\C ^2 -L) \cong \pi _1(\C ^2 -\cL')\oplus \Z$$

Fan \cite{Fa2} has proved a similar result for the projective \fg:

\begin{prop}[{\cite[Lemma 3.1]{Fa2}}] \label{Fan_2la}
Let $\cL_1, \cL_2 \sbs \C\PP ^2$ be arrangements of $n_1,n_2$
lines respectively, such that $\cL_1$ and $\cL_2$
intersect at $n_1 n_2$ distinct points in $\C\PP ^2$. Let $L$ be a line
which does not contain any of these intersection points,
and let $\cL=\cL_1 \cup \cL_2 \cup L \sbs \C\PP^2$.
Then:
$$\pi _1 (\C\PP ^2 - \cL) \cong \pi _1 (\C\PP ^2 -(\cL_1 \cup L)) \oplus \pi _1 (\C\PP ^2 -(\cL_2 \cup L)).$$
\end{prop}

In particular, if $\cL=\cL'  \cup  L $ is a line arrangement,
where $L$ is a simple line in $\cL$, then we get the following corollary:

\begin{cor}\label{simpline}
We have that:
$$\pi _1(\C\PP ^2 - \cL) \cong \pi_1(\C\PP^2 - \cL ') \oplus \Z. $$
\end{cor}

\begin{proof}
Let $\cL ' = \cL''  \cup L_1$ where $L _1$ is any line in $\cL$. Then:
\begin{eqnarray}
\pi _1(\C\PP ^2 - \cL) & = & \pi _1(\C\PP ^2 - (\cL'' \cup L \cup L _1)) \nonumber \\
    & \stackrel{{\rm Prop.\ \ref{Fan_2la}}}{\cong} & \pi_1(\C\PP^2 - (\cL '' \cup L _1)) \oplus \pi_1(\C\PP^2 - (L \cup L _1)) \nonumber \\
    & \cong & \pi_1(\C\PP^2 - (\cL '' \cup L_1)) \oplus \Z \nonumber \\
    & = &  \pi_1(\C\PP^2 - \cL ') \oplus \Z  \nonumber
\end{eqnarray}
\end{proof}

\subsection{Arrangements with at most $6$ lines}
Fan (\cite{Fa1}, \cite{Fa2}) has classified the arrangements with $6$
lines which have more than two multiple points, and proved that in these cases
the projective fundamental group is determined by the signature
of the arrangement.

Together with Theorem \ref{Fan}, he proves that the projective
fundamental group of a line arrangement with at most $6$ lines
depends only on its signature \cite{Fa2}. This classification is
repeated here at Subsection \ref{6-large}.

\section{Actions on wiring diagrams} \label{actions}

Fix a signature $S=[2^{n_2}3^{n_3}...]$. Denote by $W_S$ the set of all lists of
Lefschetz pairs with that given signature, for which the associated
wiring diagram has the unique intersection property (see Subsection \ref{signature}).
Note that the number of lines and intersection points is determined by $S$
(there are $p=\sum n_k$ points, and the number of lines $\ell$ is determined by
Equation (\ref{SUIP})).

In this section, we shortly describe two equivalence relations 
and three actions on $W_S$ or its quotient sets which were
introduced in \cite{GTV}. The actions are motivated by the
natural isometries of the plane: reflection (denoted by $\tau$),
rotation ($\mu$) and `shift of infinity' ($\si$). In \cite{GTV}
it is shown that these equivalence relations and actions preserve
the incidence lattice and the affine and projective fundamental
groups.

\smallskip

We say that two Lefschetz pairs
$\Lpair{a}{b}$ and $\Lpair{c}{d}$ are \defin{disjoint } if the integral
 segments $[a,b]$ and $[c,d]$
are disjoint.

Two lists of Lefschetz pairs ${\bf a, b} \in W_S$ are
\defin{equivalent}, denoted by ${\bf a} \equiv \bb$,
if it is possible to reach from ${\bf a}$ to $\bb$ by switching
adjacent disjoint pairs. This is easily seen to be an equivalence
relation on $W_S$.

Next, we define the relation $\arel$. 
In order to simplify notation, we write $c+\Lpair{a}{b}$ for the
\LP\ $\Lpair{c+a}{c+b}$. Also, if $\ba$ is a list of \LPs, by
$c+\ba$ we mean adding $c$ to all the pairs in $\ba$. Suppose $1
\leq c \leq c+i \leq c+t \leq \ell$.
We define the following two lists of \LPs:
\begin{eqnarray*}
\tru[]{i}{t} & = & (\Lpair{i}{i+1}, \cdots, \Lpair{t-1}{t},\Lpair{0}{t-1},\\
& & \Lpair{t-1}{t}, \cdots, \Lpair{t-i}{t-i+1}),\\
\trd[]{i}{t} & = & (\Lpair{i-1}{i},\cdots,\Lpair{0}{1},
\Lpair{1}{t},\\
& & \Lpair{0}{1}, \cdots,\Lpair{t-1-i}{t-i}).
\end{eqnarray*}

Now, we say that ${\bf a} \arel {\bf b}$, if
there is a chain of replacements of
$\tru{i}{t}$ by $\trd{i}{t}$, or \viceversa, which goes
from $\ba$ to $\bb$.

We move to define the first action, $\tau$. If a line arrangement
$\cL$ induces a list
$\ba$ of Lefschetz pairs, then reflection with respect to a line
perpendicular to the guiding line will invert the order of the
pairs in $\ba$. More generally, we define $\tau$ as follows:

\begin{defn}
For ${\bf a} = (\Lpair{a_1}{b_1},\Lpair{a_2}{b_2}, \cdots , \Lpair{a_p}{b_p}) \in W_S$, let
$$\tau({\bf a}) := (\Lpair{a_p}{b_p},\Lpair{a_{p-1}}{b_{p-1}}, \cdots , \Lpair{a_1}{b_1}).$$
\end{defn}

Note that if ${\bf a} \equiv {\bf b}$, then $\tau({\bf a}) \equiv \tau({\bf b})$.

The next action, $\mu$, is motivated by the behavior of the induced
wiring diagram when rotating the line arrangement.
We decompose the list of pairs {\bf a} into three disjoint
sublists $L_+,L_-$ and $L_0$. For this, we view every $\Lpair{a_i}{b_i}$
as a permutation acting on the indices
$\{ 1, \cdots, \ell \}$, where  $\Lpair{a_i}{b_i}$ sends every $a_i \leq t \leq b_i$ to
$a_i+b_i-t$, and leaves the other indices fixed.

Now, set $x=1$ and $L_+=L_-=L_0=\emptyset$. For each $i=1,\dots ,p$, act as follows.
If $a_i \leq x \leq b_i$, add $\Lpair{a_i}{b_i}$ to the list $L_0$, and set $x=a_i+b_i-x$.
If $x>b_i$, add $\Lpair{a_i}{b_i}$ to the list $L_-$, and if $x<a_i$, add $\Lpair{a_i}{b_i}$
to the list $L_+$. Continue to the next value of $i$.

During this procedure, the line numbered $1$ at $\infty$ always carries
the local index $x$. As a result, whenever $a_i \leq x \leq b_i$, we actually
have $x=a_i$ (since otherwise the line with local number $a_i$ would have to intersect
the first line at some point before, and again at point number $i$, contradicting the unique intersection property).

If follows that the intersection points which were combined as $L_0$ all lay
on the first line, the points in the list $L_+$ are above that line,
and the ones in $L_-$ are below it.

It is now obvious that we can reorder ${\bf a}$ into an equivalent list ${\bf a}'$,
which can be written as an {\it ordered} union ${\bf a}'= L_+ \cup L_0 \cup L_-$.
Thus we can define an action $\mu$ as follows.

\begin{defn}\label{def_mu}
Let ${\bf a} \in W_S$ be a list of Lefschetz pairs,
$${\bf a} = (\Lpair{a_1}{b_1},\Lpair{a_2}{b_2}, \cdots , \Lpair{a_n}{b_n}).$$
Decompose ${\bf a} \equiv L_+ \cup L_0 \cup L_-$ by the procedure described above,
so that we can write:
$$L_+ = (\Lpair{a_1}{b_1}, \cdots , \Lpair{a_u}{b_u}),$$
$$L_0 = (\Lpair{a_{u+1}}{b_{u+1}}, \cdots , \Lpair{a_v}{b_v}),$$
$$L_- = (\Lpair{a_{v+1}}{b_{v+1}}, \cdots , \Lpair{a_n}{b_n}).$$
where $1 \leq u < v \leq n$.

We subtract one from the indices in $L_+$, invert $L_0$ and add one to the indices in $L_-$; $\mu({\bf a})$ is defined by:
\begin{eqnarray*}
\mu({\bf a}) & := & (\Lpair{a_1-1}{b_1-1}, \cdots , \Lpair{a_u-1}{b_u-1}, \\
& & \Lpair{a_v}{b_v},\Lpair{a_{v-1}}{b_{v-1}}, \cdots ,\Lpair{a_{u+1}}{b_{u+1}}, \\
&  & \Lpair{a_{v+1}+1}{b_{v+1}+1}, \cdots , \Lpair{a_n+1}{b_n+1}).
\end{eqnarray*}
\end{defn}

Since the decomposition ${\bf a} \equiv L_+ \cup L_0 \cup L_-$ is unique,
it is obvious that for ${\bf a} \equiv {\bf b}$ we get $\mu({\bf a})= \mu ({\bf b})$.

Finally, denote by
$J$ the inversion permutation:
$J(i)=\ell+1-i, \quad 1 \leq i \leq \ell$.

\begin{defn}
Let $\ba = (\Lpair{a_1}{b_1},\Lpair{a_2}{b_2}, \cdots , \Lpair{a_p}{b_p}) \in W_S$.
Then:
$$\si(\ba) := (\Lpair{a_2}{b_2},\Lpair{a_3}{b_3}, \cdots , \Lpair{a_p}{b_p},\Lpair{J(b_1)}{J(a_1)}).$$
\end{defn}

Note that $\sigma$ does not preserve equivalence
classes of $\equiv$.

\medskip

We say that two wiring diagrams ${\bf a}$ and ${\bf b}$ are
\defin{similar} if there is a chain of wiring diagrams from ${\bf
a}$ to ${\bf b}$, such that in each step we move to an equivalent
wiring diagram (with respect to $\equiv$ or $\arel$), or act with
$\tau,\mu$ or $\si$. Obviously, this is again an equivalence relation.

The following theorem
 is proven in \cite{GTV}:
\begin{thm}
Similar wiring diagrams have the same incidence lattice, and the
same affine and projective fundamental groups.
\end{thm}

In the next section we describe an algorithm used to enumerate and find
representatives of all the similarity classes for a given signature,
for all the signatures of up to $8$ points. The detailed  results are then
given in sections \ref{6lines}--\ref{8lines}.

\section{Classification of similarity classes} \label{alg}

As explained in previous sections, our aim is to classify the wiring diagrams
of up to $8$ wires (which satisfy the unique intersection property)
up to similarity (which is the transitive closure of $\equiv, \arel$
and the relations induced by the actions $\tau$,$\mu$ and $\si$).
This classification is accomplished along the following lines.

\begin{enumerate}
\item First, we find the possible signatures for a line arrangement
with a given number of lines. For this we only need to solve Equation (\ref{SUIP})
for integral $n_2,n_3, \cdots$, and this is easy when $\ell$ is small. We use
Lemma \ref{largepoint} below to rule out some impossible signatures.
\item For every fixed signature, we produce a list of one representative from each
similarity class. Details on this crucial step are given in Subsection \ref{enumss}.
\item Next, we find the incidence lattice of every similarity class.
It is common to discover that the few dozens of similarity classes have very few
lattice structures.
The mechanism used for this part is described in Subsection \ref{latss}.
\item In the final step, we check that every two similarity classes which have
isomorphic lattices, also have the same affine and projective fundamental groups.
This test was performed using the software package \software{testisom}.
Details are given in Subsections \ref{comp_fund_group} and \ref{testisomss}.
\end{enumerate}

The reader will observe that our procedure classifies wiring diagrams,
and not only line arrangements. However, by a result of Goodman and Pollack \cite{GP},
the minimal wiring diagram with the unique intersection property which is not induced
by a line arrangement has $9$ wires.

\subsection{Signatures of line arrangements}

The following lemma shows that many solutions to
$\sum _{k\geq 2} n_k{{k} \choose {2}} = {{\ell} \choose {2}}$ do not occur as signatures
of a line arrangement:

\begin{lem}\label{largepoint}
Let $[2^{n_2} 3^{n_3} \cdots ]$ be a signature of a wiring diagram with $\ell$ pseudolines,
and assume there is a point in which $\ell-c$ lines meet. Then:
$$\sum  n_k {{k-1} \choose {2}} \leq {{c} \choose {2}} + {{\ell-c} \choose {2}}$$
\end{lem}
\begin{proof}
Let $P$ be the point of multiplicity $\ell-c$, and consider the remaining $c$ lines.
If they form a diagram with signature $[2^{m_2} 3^{m_3} \cdots ]$
(so that $\sum  m_k{{k} \choose {2}} = {{c} \choose {2}}$), then in the original diagram the multiplicity
of every point cannot be larger by more than $1$. Since we now count the point $P$ too,
we have that $\sum  n_k {{k-1} \choose {2}} \leq {{c} \choose {2}} + {{\ell-c} \choose {2}}$.
\end{proof}

If we do not count the point in which $\ell -c$ lines meet, then  for $c=1$,
the lemma gives $n_3+\cdots \leq 0$, so the signature is
$[2^{\ell-1} (\ell -1)^1]$. For $c=2,3,4$, the inequalities are $n_3+ \cdots \leq 1$,
$n_3+3n_4 +\cdots \leq 3$ and $n_3+3n_4 +6n_5 + \cdots \leq 6$, respectively.

\subsection{Computer enumeration of similarity classes} \label{enumss}

Fix a signature $S$ satisfying Equation (\ref{SUIP}), for example
$[2^{16} 3^{4}]$. 
(this signature has one of the largest spaces $|\Wmodequiv|$ among the signatures for which this size was computed).

The number of lines and
intersection points is determined by $S$; in the example, there
are $\ell=8$ lines and $p=20$ intersection points. We want to
list all the wiring diagrams with signature $S$, up to the action
of $\tau,\mu,\si$ and 
the equivalence relations $\equiv,
\arel$. We actually intend to generate a list, in which every
element is minimal in its similarity class (with respect to the
lexicographical order).

Recall that to a wiring diagram $\ba$, we associate a list
$(\Lpair{a_1}{b_1},\cdots,\Lpair{a_p}{b_p})$ of Lefschetz pairs, where
the multiset of differences $b_i-a_i$ gives the signature.
Viewing each pair $\Lpair{a}{b}$ as the permutation sending $a
\leq x \leq b$ to $a+b-x$, the pseudoline numerated $1 \leq i
\leq \ell$ at $+ \infty$, becomes $\Lpair{a_p}{b_p} \cdots
\Lpair{a_1}{b_1} (i)$ at
$- \infty$. By Lemma \ref{UIP}, $\ba$ has the unique intersection property if and only if
$\Lpair{a_p}{b_p} \cdots \Lpair{a_1}{b_1} =J$, where $J= \Lpair{1}{\ell}$
sends $x\mapsto \ell+1-x$.

A naive approach to finding all the possible wiring diagrams would be to try every possible
lists of Lefschetz pairs $(\Lpair{a_1}{b_1},\cdots,\Lpair{a_p}{b_p})$, and  keep
only those lists for which the product is $J$. This is not practical, as there are
too many products to check (in our example $S=[2^{16}3^4]$, there are
$\binom{20}{4} \cdot 7^{16} \cdot 6^{4} = 10^{20.3}$ such products).
Instead, we use a method known as
'meet in the middle' (see e.g. \cite{MH}). This method is used to
solve the following type of problems:
$A_1,\dots,A_p$ are given (small, say $|A_i| = 5$ or $|A_i| = 7$) sets of permutations,
and we want to list all the solutions to $\pi_p \pi_{p-1} \dots
\pi_1 = J$, $\pi_i \in A_i$, for some fixed permutation $J$. The
method trades time for memory in the following way. Instead of
going through all the $|A_1|\cdots |A_p|$ possibilities for
$p$-tuples
$(\pi_1,\cdots,\pi_p)$, we fix some $1<p_0<p$, and generate a sorted list
of all the products $\pi_{p_0} \cdots \pi_1$ ($\pi_i \in A_i$).
After this half is finished, we
construct all the products $ \pi_{p_0+1}^{-1} \cdots \pi_p^{-1} J$ and look
for each of them in the list generated before.
If $\pi_{p_0} \cdots \pi_1=\pi_{p_0+1}^{-1} \cdots \pi_p^{-1} J$, then
$(\pi_1 , \cdots, \pi_p)$ is a solution to the original problem.

If, for example, the sets $A_i$ all have the same size $c$ and we take
$p_0 = p/2$, then this procedure
takes $2 c^{p/2}$ operations and memory of $c^{p/2}$ permutations, instead of
$c^p$ operations (and no memory) in the naive way. In general, if memory of
$c^{p/2}$ is available, it would be best to take $p_0 \sim p/2$.


In our problem, we would like to use the sets
$A_i = \{ [1,c], [2,c+1], \cdots, [\ell+1-c,\ell] \}$, where $c$ is the number of lines
intersecting in the $i$th point. However, we do not know in advance how many lines
intersect in each point: the signature only determines how many points of each type
there are. To overcome this obstacle, we go through all the possible
ways to break the given signature $S$ into a sum of signatures $S = S_0 + S_1$,
where $S_0$ carries  $p_0$ intersection points, and $S_1$  the remaining $p-p_0$ points.
For each dissection $S = S_0+S_1$, we go through the different
choices of the numbers of lines meeting in every point; then the sets
$A_i$ are fixed. We continue as in the standard method, and collect the
solutions from all the possible dissections.

All this suggests the following algorithm: first use the 'meet in
the middle' method (with certain variations) to list all the
wiring diagrams with signature $S$, then somehow pick the minimal
element from each similarity class. However, this method is still
not practical due to the huge number of possible wiring diagrams.
For $S = [2^{16}3^{4}]$ there are about $10^{10.6}$ such classes
(this was not directly computed; the way we got this estimate is
explained below).

Luckily, we are interested in the list of similarity classes and not in the full list
of wiring diagrams.

Fixing $1 \leq p_0 \leq p$, we define the following equivalence relation, denoted
$\halfequiv$: given a wiring diagram with $p$ intersection points, we can uniquely write it
as a composition of two wiring diagrams on the same wires,
${\bf a} = {\bf a}_2 {\bf a}_1$, where ${\bf a}_1$ contains the first $p_0$
intersection points, and ${\bf a}_2$ the other $p-p_0$. Decompose
${\bf a} = {\bf a}_2 {\bf a}_1$ and ${\bf b} = {\bf b}_2 {\bf b}_1$. Then
${\bf a} \halfequiv {\bf b}$ if and only if  ${\bf a}_1 \equiv {\bf b}_1$ and
${\bf a}_2 \equiv {\bf b}_2$. Obviously, $\halfequiv$ is an equivalence relation
(with smaller classes than those of $\equiv$, in general).

We use $<$ to denote the lexicographical order between wiring diagrams
(\ie\ the associated lists of Lefschetz pairs). Suppose
${\bf a} = {\bf a}_2 {\bf a}_1 \in W_S$ is a wiring diagram; if ${\bf a}_2' {\bf a}_1$
is also a wiring diagram with ${\bf a}_2' \equiv {\bf a}_2$ and
${\bf a}_2 < {\bf a}_2'$,
then ${\bf a}_2 {\bf a}_1 < {\bf a}_2' {\bf a}_1$. We may thus
keep only the minimal among all $\equiv$-equivalent halves
$(\Lpair{a_1}{b_1},\cdots,\Lpair{a_{p_0}}{b_{p_0}})$ which
multiply to the same permutation. The same argument holds of
course for the second half of the wiring diagrams.

To summarize, this improved 'meet in the middle' mechanism produces a list of elements,
which are minimal representatives in their classes with respect to $\halfequiv$.
Every wiring diagram in $W_S$ is $\halfequiv$-equivalent to some element in this
list.

Minimality in both halves does not guarantee that an element is minimal in its
$\equiv$-class. Indeed, such $\ba$ is not minimal exactly when
$\Lpair{a_{p_0+1}}{b_{p_0+1}} < \Lpair{a_{p_0}}{b_{p_0}}$. We simply
erase these elements, and the result is a list $\omega_S$ of minimal representatives
for all the $\equiv$-classes. For the signature $S = [2^{16}3^{4}]$, there are
$3317776$ $\equiv$-classes (and $8289348$ $\halfequiv$-classes).

Given this list, it is possible to estimate the size of $W_S$ as follows. We randomly choose some
of the representatives, compute the size of their $\equiv$-class (the complexity is linear in the
size of the class), and use the average to
estimate the average size of all the classes. This average is then multiplied by the number of
classes.

Let $[\ba]$ denotes the $\equiv$-class of $\ba$. Recall that $\ba,{\bf b} \in \omega_S$ are
similar if it is possible to reach from an element of $[\ba]$ to an element of
$[\bb]$ by action of $\tau,\mu$ or $\si$, or by replacing a wiring diagram
${\bf c}$ by ${\bf c'} \equiv {\bf c}$ or ${\bf c'} \arel {\bf c}$.

For $\ba \in W_S$, we denote by $\ba _0$ the minimal element in $[\ba]$
(with respect to the lexicographical order). For every $\ba \in \omega_S$, we compute
$\tau (\ba)_0$ and $\mu (\ba)_0$ which are in $\omega_S$. We also compute the set
$\si ( [\ba] )_0$: $\si$ is not well-defined on $\Wmodequiv$, but for $\ba ' \equiv \ba$,
$[\si (\ba ')]$ depends only on the leftmost pair in $\ba '$. Thus, $\si ( [ \ba ] )_0$ is
of size no more than $p$ (but usually much less). Finally, we find $\bb_0$ for every
$\bb \arel \ba$. This is done by locating lists of pairs of the form
$\tru{i}{t}$ or $\trd{i}{t}$,
which are consecutive in some member of
$[\ba]$, and replacing $\tru{i}{t}$ by $\trd{i}{t}$, or \viceversa.

We generate a graph on $\omega_S$ as the set of vertices, where $\ba \in \omega_S$ is
connected to all the elements $\tau (\ba)_0, \mu(\ba)_0, \si( [\ba] )_0$
and $\bb_0$ for $\bb \arel \ba '$, $\ba ' \in [\ba]$.

It remains only to list the minimal representatives for the
connected components of this graph by standard depth-first
search. The result is the required list of minimal
representatives of the similarity classes.

The running time of this algorithm depends of course on the size of $W_S$.
The worst case for arrangements of up to $8$ lines was $[2^{19} 3^3]$,
which took several CPU hours on a medium-size Unix machine;
the size of $W_S$ is about $10^{11.75}$, 
and
there were $5$ similarity classes. 

To demonstrate the relative influence of the relation $\arel$ and
the actions $\mu,\tau,\si$, we give below the number of classes
of $\Wmodequiv$ (for $S = [2^{16}3^4]$) under the various sets of
relations. For example, the forth line tells us that under the
transitive closure of $\mu,\arel$ and $\equiv$, there are $10948$
classes.

\begin{center}
\begin{tabular}{|cccc||c|}
            \hline
             $\sigma$
        & $\tau$
        & $\mu$
        & $\arel$
        & Number of classes \\
            \hline
- & - & - & - & 3317776 \\
- & - & - & + & 306328 \\
- & - & + & - & 207361 \\
- & - & + & + & 10948  \\
- & + & - & - & 1658965  \\
- & + & - & + & 125186  \\
- & + & + & - & 103719 \\
- & + & + & + & 5543  \\
+ & - & - & - & 35100 \\
+ & - & - & + & 460 \\
+ & - & + & - & 723 \\
+ & - & + & + & 18  \\
+ & + & - & - & 17644 \\
+ & + & - & + & 259 \\
+ & + & + & - & 723 \\
+ & + & + & + & 18 \\
\hline
\end{tabular}
\end{center}

The redundancy of $\tau$ in the presence
 of $\sigma$ and $\mu$ is
explained by \cite[Cor. 5.3]{GTV}.

\subsection{Lattice comparison}\label{latss}

Let $\ba =(\Lpair{a_1}{b_1}, \cdots, \Lpair{a_p}{b_p})$ be a list of Lefschetz pairs.
The incidence lattice of $\ba$ is determined by the sets of intersection points
on each of the wires. To compute the lattice, we compute the global numeration
of the wires at each point, by applying the permutation $\Lpair{a_i}{b_i} \cdots \Lpair{a_1}{b_1}$
on the list $1, \cdots, \ell$ for every point $i$. In what follows,
we view the incidence lattice as the bipartite graph whose vertices are the lines and points.

The problem of checking if two bipartite graphs are isomorphic is
known to be difficult in general (see \cite[p. 285]{GJ}). However
(as is often the case), one can get quite good results by
heuristic methods.

We hold a graph as a matrix $M$ in the usual way: $M_{ij}=1$
if the point number $j$ is on line number $i$, and $M_{ij}=0$ otherwise.

The idea is to sort the rows and the columns of the matrix
alternatively. After doing it several times, we compare the two sorted matrices.
If we get equal matrices, the corresponding lattices are equivalent.

Using various sorting orders, we are able to classify the lattices into
few families. Then, we  manually check if these families are indeed distinct.

\subsection{Computation of the fundamental group}\label{comp_fund_group}

The computation of the fundamental group of the complement of real line
arrangements is done in two steps: the Moishezon-Teicher algorithm
computes the skeletons associated to the arrangement,
and the van Kampen theorem produces from the skeletons
a finite presentation of the affine fundamental group.
One more relation gives the projective fundamental group.

The theoretical background for both steps can be found in
\cite{MoTe1}, \cite{VK} and \cite{GaTe}. Here we only sketch the
ideas from an algorithmic point of view. Full details are given
in Section 3 of \cite{GTV}.

Let $\ba = (\Lpair{a_1}{b_1},\cdots,\Lpair{a_p}{b_p})$ be a list of Lefschetz pairs.
In order to compute the skeleton $s_i$ associated to the
$i$th intersection point $\Lpair{a_i}{b_i}$
, we start with the configuration in Figure \ref{fig3}, in which
segments connect the points from
$a_i$ to $b_i$.

\FIGURE{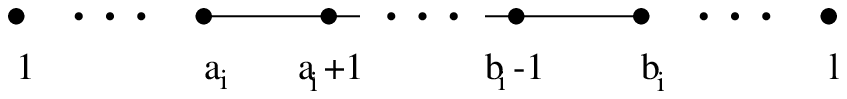}{0.7cm}{fig3}

and apply the Lefschetz pairs $\Lpair{a_{i-1}}{b_{i-1}}, \cdots,\Lpair{a_1}{b_1}$.
A Lefschetz pair  $\Lpair{a_j}{b_j}$ acts by rotating the region
from $a_j$ to $b_j$ by $180^{\circ}$ counterclockwise without affecting
any other points.

For example, consider the list
$(\Lpair{2}{3},\Lpair{2}{4},\Lpair{4}{5},\Lpair{1}{3},\Lpair{3}{4})$.
The initial skeleton for $\Lpair{3}{4}$ is given in Figure
\ref{fig4}.
\FIGUREx{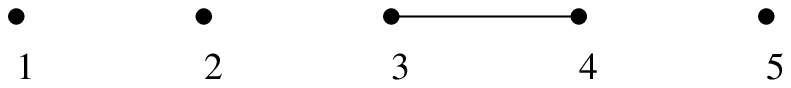}{6cm}{fig4}
Applying $\Lpair{1}{3}$ and then
$\Lpair{4}{5}$, we get the skeleton in Figure \ref{fig5}.
\FIGUREx{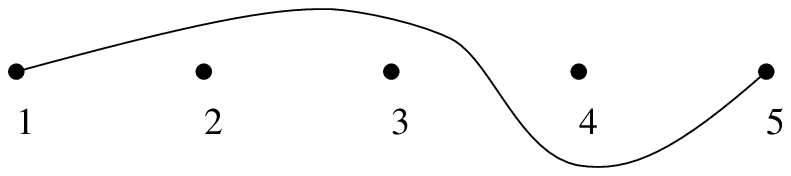}{6cm}{fig5}
Then, applying  $\Lpair{2}{4}$ yields the skeleton of Figure
\ref{fig6}, and finally, acting with  $\Lpair{2}{3}$ we get the skeleton in Figure
\ref{fig8}.
\FIGUREx{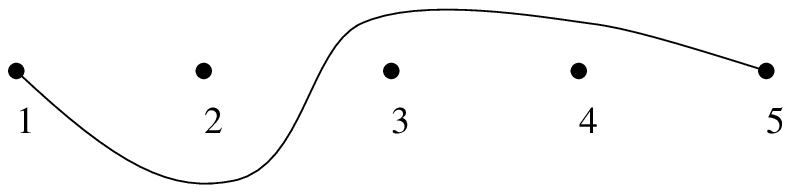}{6cm}{fig6}

\FIGUREx{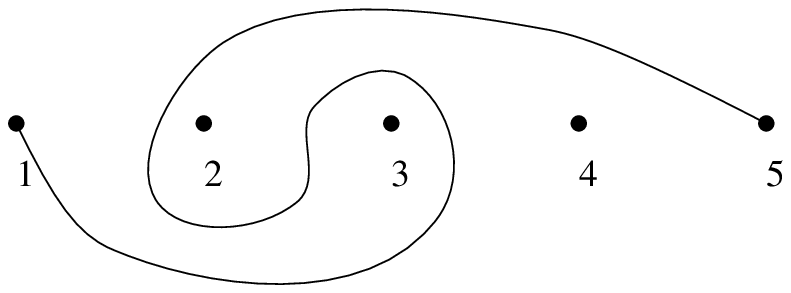}{6cm}{fig8}

Since we had to compute hundreds of fundamental groups, this
procedure had to be programmed. The algorithm was implemented as
follows: we decode a path as a list of the points near which the
path passes, with special marks for 'pass above' and 'pass
below'. For example, the path in Figure \ref{fig8} is decoded as
$1^{\circ} 2^{-} 3^{-} 3^{+} 2^{-} 2^{+} 3^{+} 4^{+} 5^{\circ}$.
Details on the implementation can be found in \cite{GKT} and \cite{Gthesis}.
\medskip

The fundamental group is generated by $\Ga_1,\dots,\Ga_\ell$, where $\Ga_j$
corresponds to loop around the $i$th point which is based on a fixed point outside the arrangement.

Assume the skeleton $s_i$ is induced by a simple intersection
point. Then the induced relation on the $\Ga_j$ is computed as
follows. We cut the skeleton just after its connection to its
start point, and denote the associated generator by $a_1$. The
second element is computed by going from the starting point
around the end point: whenever the path passes above a point we
multiply by the generator $\Ga_j$ associated to that point or its
inverse, depending on the direction of the path (we ignore the
points which we pass from below). When we reach the end point, we add
the generator associated to this point, and then go all the way
back by adding the corresponding inverse elements. This results
with an element
$a_2$, and the relation induced by $s_i$ is $a_1 a_2 = a_2 a_1$.

For example, the skeleton in Figure \ref{fig5} induces the relation
$$\Ga _3 \Ga _2 \Ga_1 \Ga _2^{-1} \Ga _3^{-1} \cdot \Ga_5 = \Ga_5 \cdot \Ga _3 \Ga _2 \Ga_1 \Ga _2^{-1} \Ga _3^{-1}.$$
Details can be found in \cite{GaTe}.

In case the skeleton $s_i$ corresponds to a multiple intersection
point, we get a set of relations of the form $a_c \dots a_2 a_1 =
a_1 a_c \dots a_2 = \dots = a_{c-1} \dots a_1 a_c$, where
$a_r$ is computed from the (same) starting point around the
$r$-th point, as described above.

To obtain the projective fundamental group, we add the relation
$\Ga_\ell \cdots \Ga_2 \Ga_1 = 1$.

\subsection{Checking isomorphisms between fundamental groups}\label{testisomss}
The problem of testing if two groups (presented in terms
of generators and relations) are isomorphic is undecidable
(see \cite{Rotman}).
However, the fact that the problem is undecidable in general
does not prevent one from trying to solve it in concrete cases.
Indeed, there is a software package called \software{testisom},
written by D.~Holt and S.~Rees,
which attempts to prove or disprove that two given groups are isomorphic.
It tries to guess an isomorphism
using the Knuth-Bendix algorithm (which generates part of Cayley graph
of the groups), or prove the groups are non-isomorphic by computation
of Abelian sections,
mainly of the lower central series, or counting projections onto
certain small groups.

For full description of this package, see \cite{EHR},\cite{HR}.

\section{Arrangements with up to $6$ lines} \label{6lines}

In this section we classify the arrangements with up to $6$ lines,
up to lattice isomorphism.

\subsection{Possible signatures for $\ell \leq 6$ lines}

The following proposition lists the different possible signatures
for arrangements with up to $6$ lines.

\begin{prop}\label{signatures-6}
The possible signatures for arrangements with up to $6$ lines (sorted by the number
of intersection points) are:
\begin{itemize}
\item[(1)]  $[2^1]$ (for $2$ lines),
\item[(2)]  $[3^1]$, $[2^3]$ (for $3$ lines),
\item[(3)]  $[4^1]$, $[2^3 3^1]$, $[2^6]$ (for $4$ lines),
\item[(4)]  $[5^1]$, $[2^4 4^1]$, $[2^4 3^2]$, $[2^7 3^1]$, $[2^{10}]$ (for $5$ lines), and
\item[(5)]  $[6^1]$, $[2^5 5^1]$, $[2^3 3^4]$,
$[2^6 3^1 4^1]$, $[2^6 3^3]$, $[2^9 4^1]$, $[2^9 3^2]$, $[2^{12} 3^1]$ and $[2^{15}]$ (for $6$ lines).
\end{itemize}
\end{prop}

\begin{proof}
There are some solutions to
$\sum _{k\geq 2} n_k{{k} \choose {2}} = {{\ell} \choose {2}}$
except for the ones listed above. For $\ell=4$ we have $[3^2]$,
which is forbidden by the case $c=1$ of Lemma
\ref{largepoint}. For
$\ell=5$, we get $[2^1 3^1 4^1]$ and $[2^1 3^3]$. The first
is forbidden by Lemma \ref{largepoint}, case $c=1$, and the
second by Lemma \ref{largepoint}, case $c=2$.

For $\ell = 6$, we get $[2^2 3^1 5^1]$, $[3^1 4^2]$, $[2^3 4^2]$,
$[3^3 4^1]$, $[2^3 3^2 4^1]$ and $[3^5]$. The first five options
are forbidden by Lemma \ref{largepoint}. As for $[3^5]$, this is
impossible since then a line $L$ would meet five other lines, two
in each intersection point; but $2$ does not divide $5$.
\end{proof}

In order to present the results of the classification, we use
diagrams based on the graphs of multiple points. In these
diagrams we omit all the simple points and the lines on which
there is at most one multiple points (these are easily put back
in place if one wants to reconstruct the arrangements from our
diagrams). We denote a multiple point of multiplicity
$3$ by a black dot, and a multiple point of higher multiplicity
by a bigger black dot with a number indicating the multiplicity.
If more than two points belong to the same line, we draw the line
thicker than the usual. From this graph, one can recover the
lattice of the arrangement, and construct explicitly the
corresponding line arrangement. For example, the arrangements of
Figure \ref{fig2} are represented by the graphs given in Figure
\ref{fig2rep}.
\FIGUREx{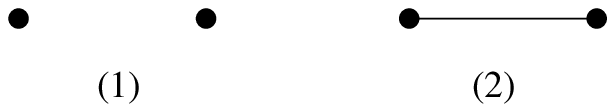}{5cm}{fig2rep}

If there are at most two multiple points, then the affine and projective
fundamental groups can be computed for every lattice by Theorems \ref{Fan-} and \ref{Fan}. These groups are given in Subsection \ref{6-small}.
For the other signatures we do not give the groups, even though their presentations were
computed for all the similarity classes. The classification of this case is
given in Subsection \ref{6-large}. We keep the same structure in Sections
\ref{7lines} and \ref{8lines}.

\subsection{Signatures with at most two multiple points}\label{6-small}

The following table summarizes the possible lattices for signatures in which
there are at most two multiple points. We give only the projective fundamental
group, and note that by Theorems \ref{Fan} and \ref{Fan-}, in these cases
the affine fundamental group is given by
\begin{equation}\label{projaff}
\pi _1 (\C^2 - \cL) \cong \pi _1 (\C\PP ^2 - \cL) \oplus \Z.
\end{equation}

\begin{center}
\begin{tabular}{|c|c|c|c|}
            \hline
             $\ell$ &Signature & Possible\ lattices &    Projective\\
                     &  &                           &  fundamental group\\
            \hline
            \hline
$2$ & $[2^1]$   &  \Empty & $\Z$ \\
            \hline
$3$ & $[3^1]$   & \tabfig{0.15cm}{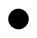} & $ \F _2$ \\
            \hline
$3$ & $[2^3]$ &  \Empty & $\Z^2$ \\
            \hline
$4$ & $[4^1]$   & \tabfig{0.5cm}{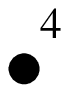} & $\F _3$ \\
            \hline
$4$ & $[2^3 3^1]$ & \tabfig{0.15cm}{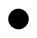} &
                                   $\Z \oplus \F _2$ \\
            \hline
$4$ & $[2^6]$ & \Empty & $\Z ^3$ \\
            \hline
$5$ & $[5^1]$  & \tabfig{0.5cm}{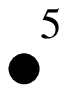} & $\F _4$ \\
            \hline
$5$ & $[2^4 4^1]$ & \tabfig{0.5cm}{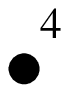} & $\Z \oplus \F _3$ \\
            \hline
$5$ & $[2^4 3^2]$  & \tabfig{0.15cm}{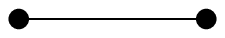}&
                       $\F _2 \oplus \F _2$ \\
            \hline
$5$ & $[2^7 3^1]$  & \tabfig{0.15cm}{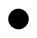} & $\Z ^2 \oplus \F _2$ \\
            \hline
$5$ & $[2^{10}]$  & \Empty & $\Z ^4$ \\
            \hline
$6$ & $[6^1]$       &  \tabfig{0.5cm}{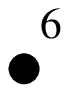} & $\F _5$ \\
            \hline
$6$ & $[2^5 5^1]$   & \tabfig{0.5cm}{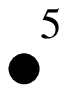} & $\Z \oplus \F _4$ \\
            \hline
$6$ & $[2^6 3^1 4^1]$ & \tabfig{0.5cm}{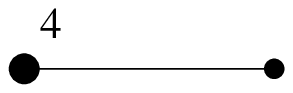} & $\F _2 \oplus \F _3$ \\
            \hline
$6$ & $[2^9 4^1]$   & \tabfig{0.5cm}{lat_4_1.ps} & $\Z^2 \oplus \F _3$ \\
            \hline
$6$ & $[2^9 3^2]$ & \tabfig{0.7cm}{lat_6_12.ps} &
                                   $\Z \oplus \F _2 \oplus \F _2$ \\
            \hline
$6$ & $[2^{12} 3^1]$ & \tabfig{0.15cm}{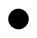} & $\Z ^3 \oplus \F _2$ \\
            \hline
$6$ & $[2^{15}]$  & \Empty & $\Z ^5$ \\
            \hline
\end{tabular}
\end{center}

\subsection{Signatures with more than two multiple points}\label{6-large}

In the following table, we summarize the possible lattices for
signatures with more than two multiple points:

\begin{center}
\begin{tabular}{|c|c|c|c|}
            \hline
             Signature & Possible\ lattices & $|\Wmodequiv|$ & Similarity\ classes\\
            \hline
            \hline
$[2^3 3^4]$ & \newline \ \tabfig{1cm}{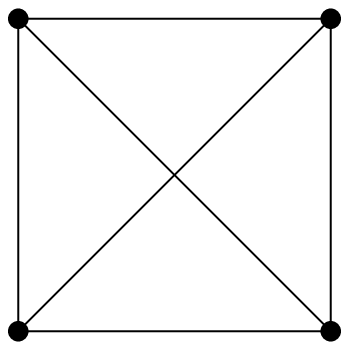}  & 16 & 1\\
            \hline
$[2^6 3^3]$ & \tabfig{1cm}{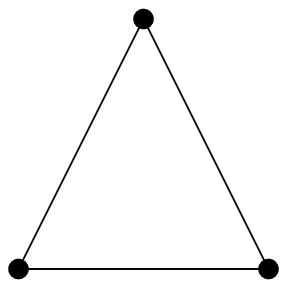}  & 304 & 2\\
            \hline
\end{tabular}
\end{center}

\section{Arrangements with 7 lines} \label{7lines}

In this section we classify the arrangements with $7$ lines up to
lattice isomorphism.

\subsection{Possible signature for $\ell=7$}

\begin{prop}\label{signatures-7}
The possible signatures for arrangements with $7$ lines are:
$[7^1]$, $[2^6 6^1]$, $[2^3 3^6]$, $[2^8 3^1 5^1]$, $[2^6 3^3 4^1]$, $[2^9 4^2]$,
$[2^6 3^5]$, $[2^9 3^2 4^1]$, $[2^{11} 5^1]$, $[2^9 3^4]$, $[2^{12} 3^1 4^1]$,
$[2^{12} 3^3]$, $[2^{15} 4^1]$, $[2^{15} 3^2]$, $[2^{18} 3^1]$, $[2^{21}]$.
\end{prop}

\begin{proof}
Except for the ones listed above, the only solutions to
$\sum _{k\geq 2} n_k{{k} \choose {2}} = {{7} \choose {2}}$ are:
$[4^1 6^1]$, $[3^2 6^1]$, $[2^3 3^1 6^1]$ which are forbidden by the case $c=1$ of Lemma
\ref{largepoint}; $[2^1 5^2]$, $[2^2 3^1 4^1 5^1]$, $[2^5 4^1 5^1]$, $[2^2 3^3 5^1]$,
$[2^5 3^2 5^1]$ which are forbidden by the case $c=2$ of this lemma;
$[3^1 4^3]$, $[2^3 4^3]$, $[3^3 4^2]$, $[2^3 3^2 4^2]$, $[2^6 3^1 4^2]$, $[3^5 4^1]$,
$[2^3 3^4 4^1]$ and $[3^7]$.

The first seven options are excluded by the case $c=3$ of the lemma,
since we must have $n_3+3n_4 \leq 3$ (note that here $n_4$ counts the points of
multiplicity $4$ beyond the first one). As for $[3^7]$, this is the projective
geometry on $7$ points, which is well known not to be embedded 
in the affine plane.
\end{proof}

\medskip

We now describe the different line arrangements in details. We
use the encoding of lattices into diagrams, described in section
\ref{6lines}.

\subsection{Signatures with at most two multiple points}

The following table lists the possible lattices for the signatures
of arrangements with $7$ lines, if there are at most two multiple points.
We give the corresponding projective fundamental groups, which are computed
using Theorem \ref{Fan}.
From this one can compute the affine groups by Equation
 (\ref{projaff}).

\medskip

\begin{center}
\begin{tabular}{|c|c|c|}
            \hline
             Signature & Possible\ lattices &         Projective\\
                       &                    & fundamental group \\
            \hline
            \hline
$[7^1]$       &  \tabfig{0.5cm}{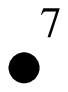} & $\F _6$ \\
            \hline
$[2^6 6^1]$   & \tabfig{0.5cm}{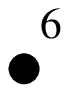} & $\Z \oplus \F _5$ \\
            \hline
$[2^8 3^1 5^1]$ & \tabfig{0.5cm}{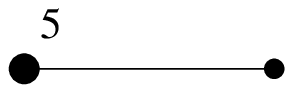} & $\F _2 \oplus \F _4$ \\
            \hline
$[2^9 4^2]$   & \tabfig{0.5cm}{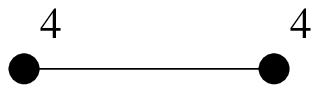} & $\F _3 \oplus \F _3$\\
            \hline
$[2^{11} 5^1]$ & \tabfig{0.5cm}{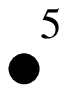} & $\Z ^2 \oplus \F _4$ \\
            \hline
$[2^{12} 3^1 4^1]$ & \tabfig{0.8cm}{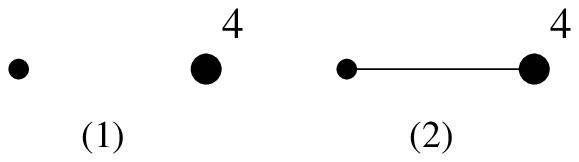} &
                                   $\Z  \oplus \F _2 \oplus \F _3$ \\
            \hline
$[2^{15} 4^1]$   & \tabfig{0.5cm}{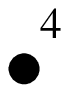} & $\Z^3 \oplus \F _3$ \\
            \hline
$[2^{15} 3^2]$ & \tabfig{0.7cm}{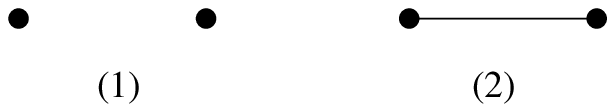} &
                                   $\Z ^2 \oplus \F _2 \oplus \F _2$ \\
            \hline
$[2^{18} 3^1]$ & \tabfig{0.15cm}{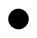} & $\Z ^4 \oplus \F _2$ \\
            \hline
$[2^{21}]$  & \Empty & $\Z ^6$ \\
            \hline
\end{tabular}
\end{center}

\subsection{Signatures with more than two multiple points}

In the following table, we summarize the possible lattices for
signatures with more than two multiple points.

\begin{center}
\begin{tabular}{|c|c|c|c|}
  \hline
Signature & Possible \ lattices & $|\Wmodequiv|$ & Similarity\ classes\\
  \hline\hline
$[2^3 3^6]$ & \tabfig{1.5cm}{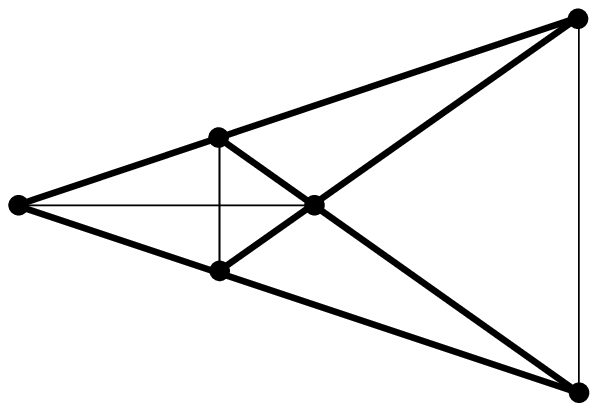} & 28 & 1 \\
\hline
$[2^6 3^3 4^1]$ & \tabfig{1.4cm}{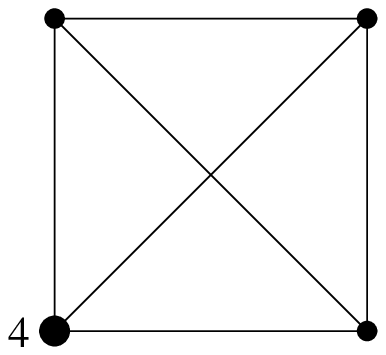} & 448 & 1 \\
\hline
$[2^6 3^5]$ & \tabfig{2.0cm}{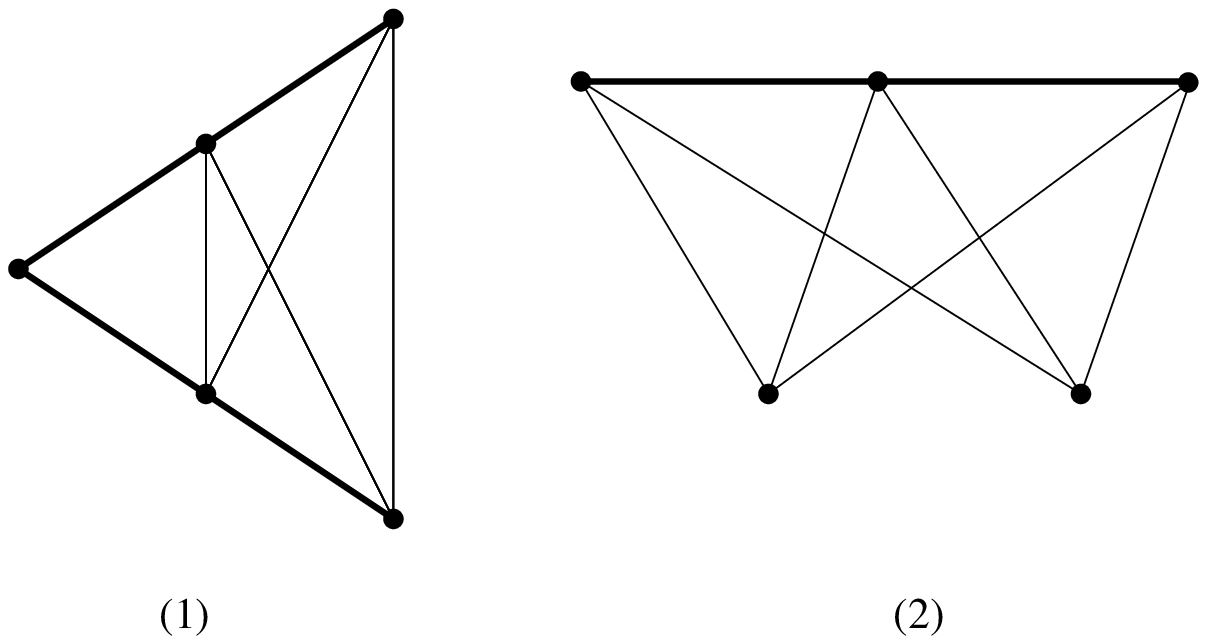} & 770 & 2  \\
\hline
$[2^9 3^2 4^1]$ & \tabfig{1.5cm}{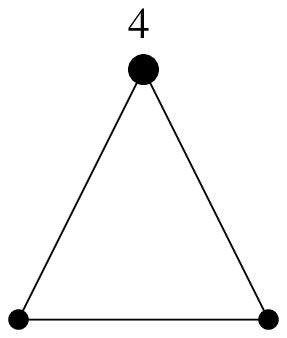} & 3864 & 3\\
\hline
$[2^9 3^4]$ & \tabfig{1.5cm}{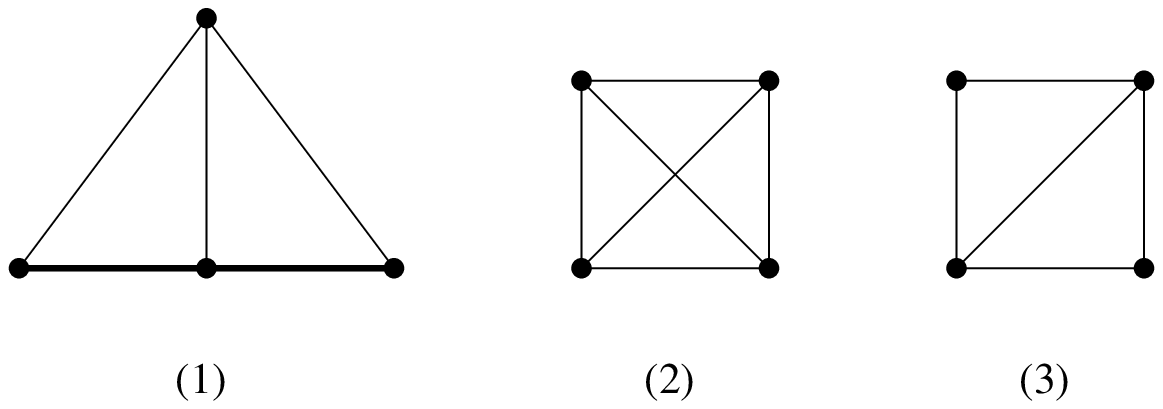} & 9534 & 6 \\
\hline
$[2^{12} 3^3]$ & \tabfig{1.5cm}{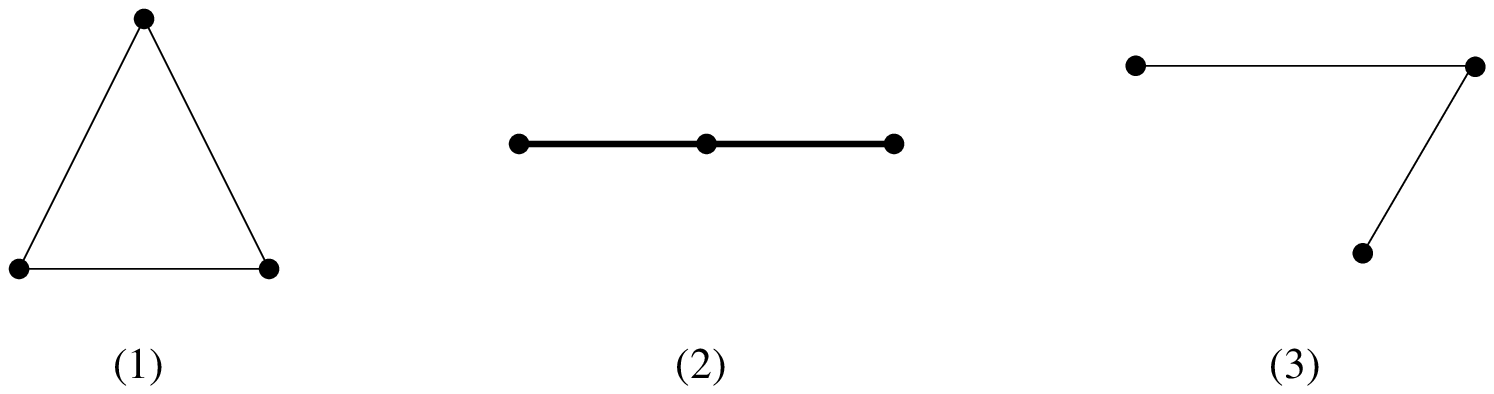} & 45458 & 4 \\
\hline
\end{tabular}
\end{center}

\section{Arrangements with 8 lines} \label{8lines}

Here we classify the arrangements with $8$ lines
up to lattice isomorphism.

\subsection{Possible signatures for $\ell = 8$}

\begin{prop}\label{signatures-8}
The following are all the possible
signatures for arrangements with $8$ lines:
$[8^1]$, $[2^7 7^1]$, $[2^4 3^6 4^1]$, $[2^{10} 3^1 6^1]$, $[2^7 3^3 4^2]$, $[2^9 3^3 5^1]$,
$[2^7 3^5 4^1]$, $[2^{12} 4^1 5^1]$, $[2^{13} 6^1]$, $[2^7 3^7]$, $[2^{10} 3^2 4^2]$,
$[2^{12} 3^2 5^1]$, $[2^{10} 3^4 4^1]$, $[2^{10} 3^6]$, $[2^{13} 3^1 4^2]$, $[2^{15} 3^1 5^1]$,
$[2^{13} 3^3 4^1]$, $[2^{13} 3^5]$, $[2^{16} 4^2]$, $[2^{18} 5^1]$, $[2^{16} 3^2 4^1]$,
$[2^{16} 3^4]$, $[2^{19} 3^1 4^1]$, $[2^{19} 3^3]$, $[2^{22} 4^1]$, $[2^{22} 3^2]$,
$[2^{25} 3^1]$, $[2^{28}]$.
\end{prop}

\begin{proof}
There are 38 solutions to the equation
$\sum _{k\geq 2} n_k{{k} \choose {2}} = {{8} \choose {2}}$ {\it except} for the above
listed ones.  Lemma \ref{largepoint} excludes three cases with $c=1$, nine cases with
$c=2$ and
fifteen cases with $c=3$. Out of the nine possibilities with a
point with multiplicity 4, the lemma (with $c=4$) excludes eight,
the remaining one is $[2^{10} 4^3]$. Suppose there are at least
two points in which four lines meet. If they do not share a
common line, the signature must be $[2^{16} 4^2]$. Otherwise,
there is one remaining line, which can meet at most 3
otherwise-simple intersection points, so we cannot have a third
point of multiplicity $4$.

The last two cases are $[2^1 3^9]$ and $[2^4 3^8]$. Assume 
signature of the form
$[2^{n_2} 3^{n_3}]$. Every line meets seven other lines, so the number of multiple points
on a line is no more than $3$. Together, there are no more than
$\frac {8 \cdot 3}{3} = 8$ multiple points  (since each point is
counted three times). This argument rules out the signature $[2^1
3^9]$. It also shows that in a diagram with signature $[2^4
3^8]$, on each line there are exactly one simple point and three
multiple points. We could not find a geometric argument to rule
this case out (see Remark \ref{8_8} below), but our algorithm
(Section
\ref{enumss}) found no wiring diagrams for this signature.
\end{proof}

\begin{rem} \label{8_8}
We would like to explain why it is difficult to show that there are no wiring diagrams
with the signature $[2^4 3^8]$. All the arguments in Propositions 
\ref{signatures-6},
\ref{signatures-7} and \ref{signatures-8} are combinatorial in nature,
actually showing that there is no arrangement of curves
with the unique intersection property which induces the forbidden signatures (with the exception
of $[3^7]$ which indeed is induced by the projective geometry on $7$ points).

However, we now show that $[2^4 3^8]$ is in fact induced by an
arrangement of curves. For that, it is enough to describe a
bipartite graph (on the sets $C$ of 'curves' and $P$ of
'points'), in which every $c_1,c_2 \in C$ 'meet' once (\ie\ are
connected to a single
$p \in P$), and such that there are $4$ points of degree (=total number of edges) $2$,
and $8$ points of degree $3$.

We let numbers denote the curves, and small Latin letters denote the points.
Consider the Cayley graph of the dihedral group $D_8$ with respect to the presentation, 
$$\sg{\al,\be,\ga \ \suchthat \ \al^2=\be^2=\ga^2=1, \ga\al\ga=\be\ga\be, \be\al\be=\al\ga\al}$$
shown in Figure \ref{cayley} (with $\be$ the generator drawn horizontally).

\FIGURE{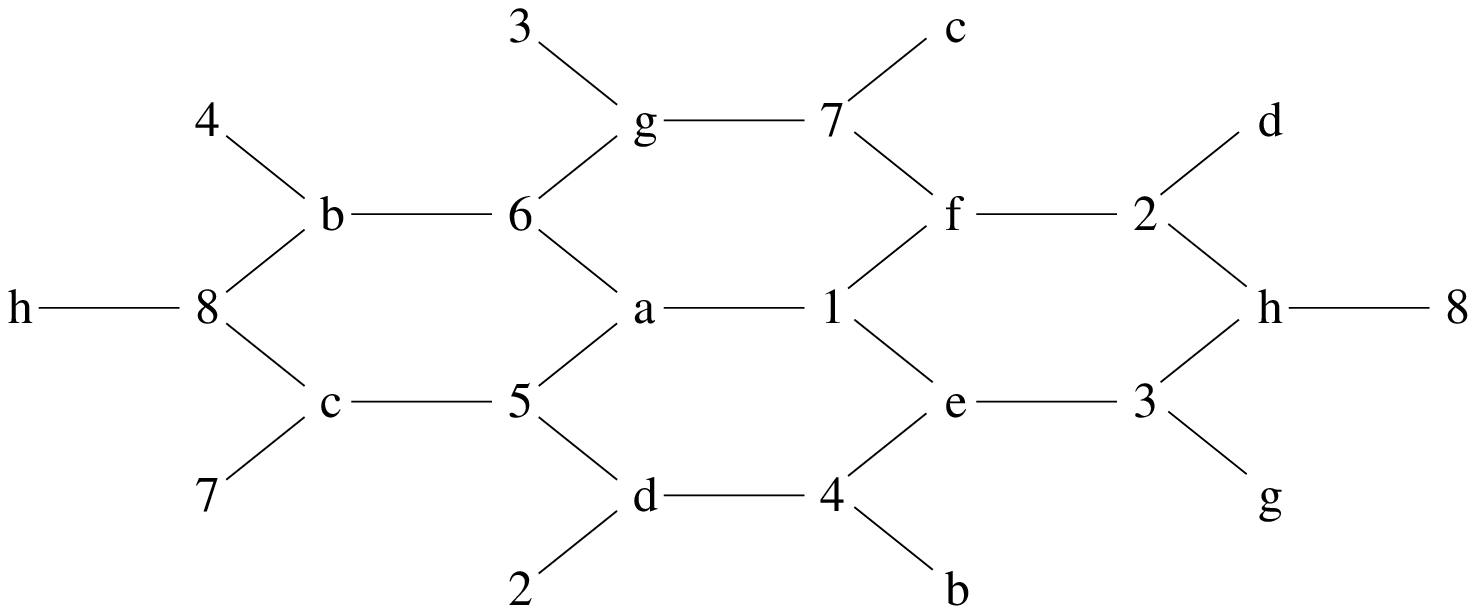}{5cm}{cayley}

This is a bipartite graph of degree $3$, on $8+8$ points
(in Coxeter's notation \cite[Figs. 11,12]{Cox}, this is $\set{8}+\set{8/3}$).
Note that the only
couples of curves which do not share a common neighbor are
$\set{1,8},\set{2,6},\set{3,5}$ and $\set{4,7}$. When we add the points connecting
these couples  we get the desired graph. This
is the only graph with these properties.
\end{rem}

\subsection{Signatures with at most two multiple points}

We now list the possible lattices and groups for signatures with
at most two multiple points for $\ell = 8$ lines. As before, we
give only the projective groups, from which one can compute the
affine groups by Equation (\ref{projaff}).

\begin{center}
\begin{tabular}{|c|c|c|}
            \hline
             Signature & Possible\ lattices &  Projective\\
                       &                    & fundamental group \\
            \hline
            \hline
$[8^1]$       &  \tabfig{0.5cm}{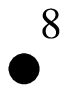} & $\F _7$ \\
            \hline
$[2^7 7^1]$   & \tabfig{0.5cm}{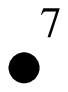} & $\Z \oplus \F _6$ \\
            \hline
$[2^{10} 3^1 6^1]$ & \tabfig{0.5cm}{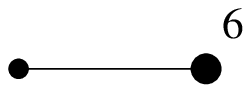} & $\F _2 \oplus \F _5$ \\
            \hline
$[2^{12} 4^1 5^1]$   & \tabfig{0.5cm}{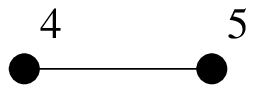} & $\F _3 \oplus \F _4$\\
            \hline
$[2^{13} 6^1]$ & \tabfig{0.5cm}{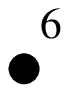} & $\Z ^2 \oplus \F _5$ \\
            \hline
$[2^{15} 3^1 5^1]$ & \tabfig{0.8cm}{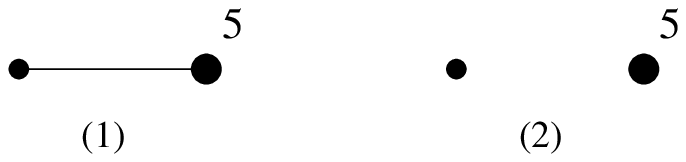} &
                                   $\Z  \oplus \F _2 \oplus \F _4$ \\
            \hline
$[2^{16} 4^2]$   & \tabfig{0.8cm}{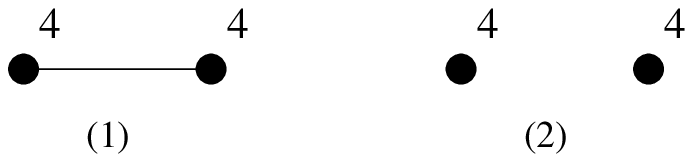} &
                                   $\Z \oplus \F _3 \oplus \F _3$ \\ 
            \hline
$[2^{18} 5^1]$ & \tabfig{0.5cm}{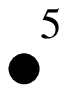} & $\Z ^3 \oplus \F _4$ \\ 
            \hline
$[2^{19} 3^1 4^1]$ & \tabfig{0.8cm}{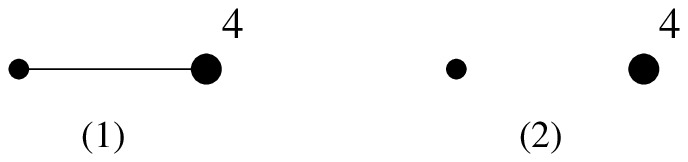} &
                                   $\Z ^2 \oplus \F _2 \oplus \F _3$ \\
            \hline
$[2^{22} 4^1]$   & \tabfig{0.5cm}{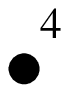} & $\Z^4 \oplus \F _3$ \\
            \hline
$[2^{22} 3^2]$ & \tabfig{0.7cm}{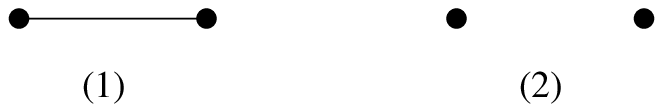} &
                                   $\Z ^3 \oplus \F _2 \oplus \F _2$ \\ 
            \hline
$[2^{25} 3^1]$ & \tabfig{0.15cm}{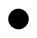} & $\Z ^5 \oplus \F _2$ \\
            \hline
$[2^{28}]$  & \Empty & $\Z ^7$ \\
            \hline
\end{tabular}
\end{center}

\subsection{Signatures with more than two multiple points}

In the following table, we summarize the possible lattices for signatures
of $\ell = 8$ lines with more than two multiple points.

\begin{center}
\begin{tabular}{|c|c|c|c|}
\hline
Signature & Possible \ lattices & $\equiv$ classes & Similarity\ classes\\
\hline\hline
$[2^4 3^6 4^1]$   & \tabfig{1.5cm}{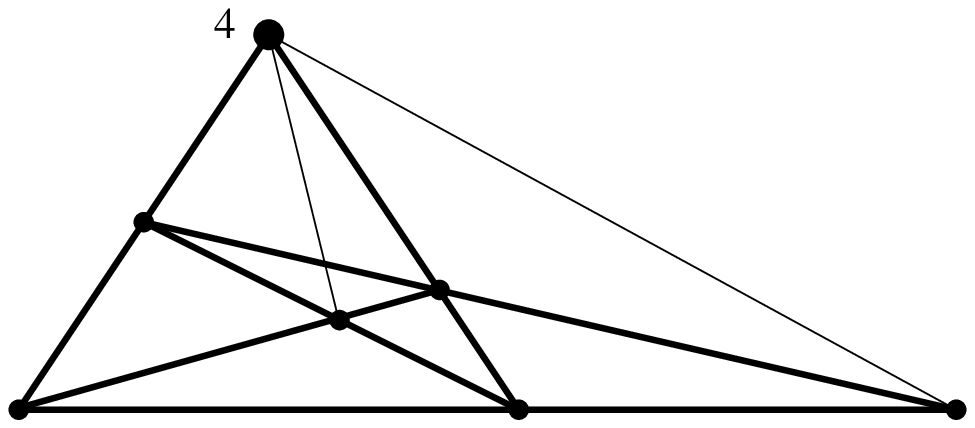} & 144 & 1 \\
\hline
$[2^7 3^3 4^2]$ & \tabfig{1cm}{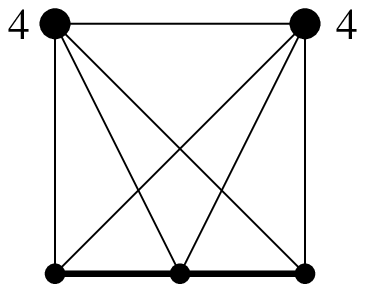} & 384 & 1 \\
\hline
$[2^9 3^3 5^1]$   & \tabfig{1cm}{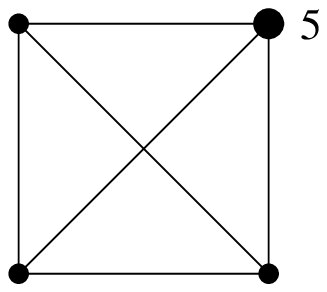} & 1792 & 2\\
\hline
$[2^7 3^5 4^1]$   & \tabfig{1.5cm}{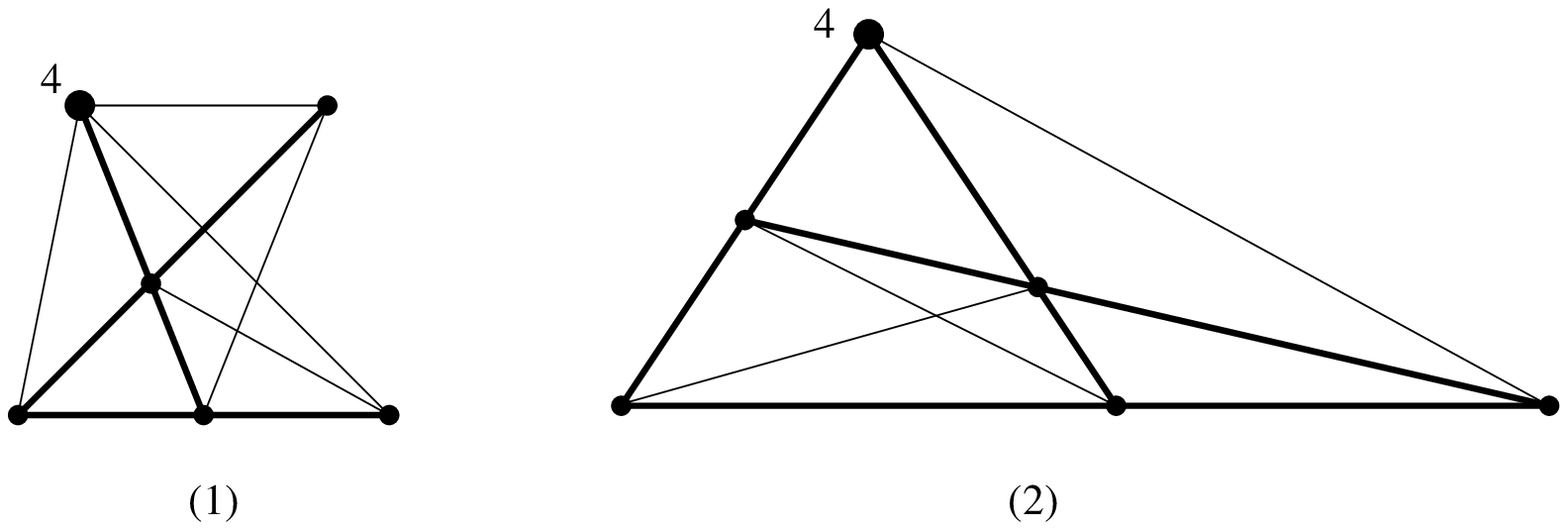} & 5024 & 4 \\ 
\hline
$[2^7 3^7]$     & \tabfig{3cm}{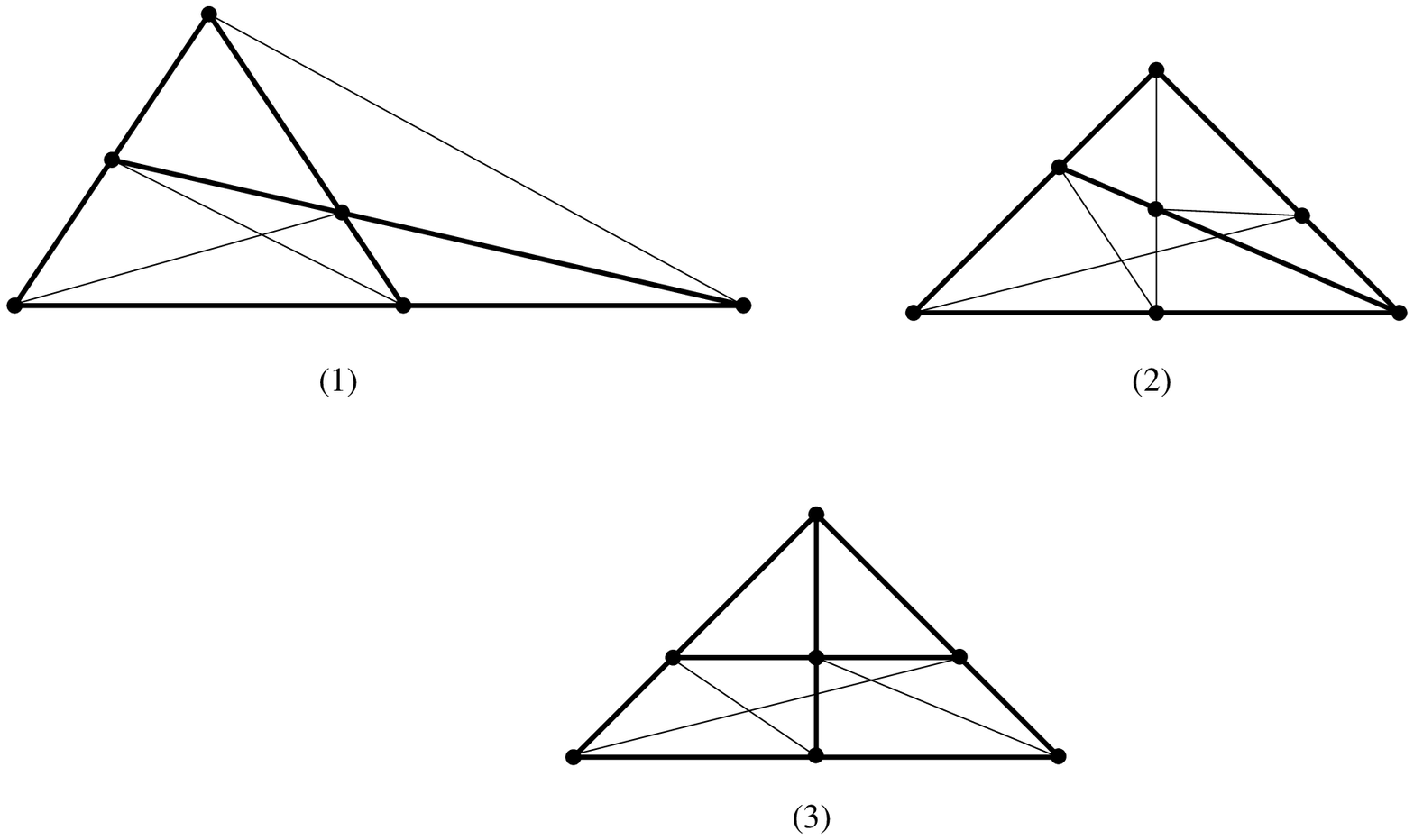} & 3920 & 4 \\ 
\hline
$[2^{10} 3^2 4^2]$ & \tabfig{1cm}{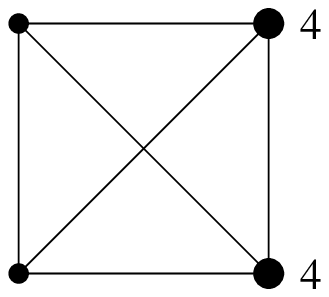} & 6528 & 4\\
\hline
$[2^{12} 3^2 5]$ & \tabfig{1cm}{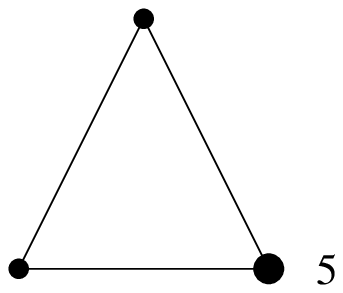} & 12096 & 4 \\
\hline
$[2^{10} 3^4 4^1]$ & \tabfig{3cm}{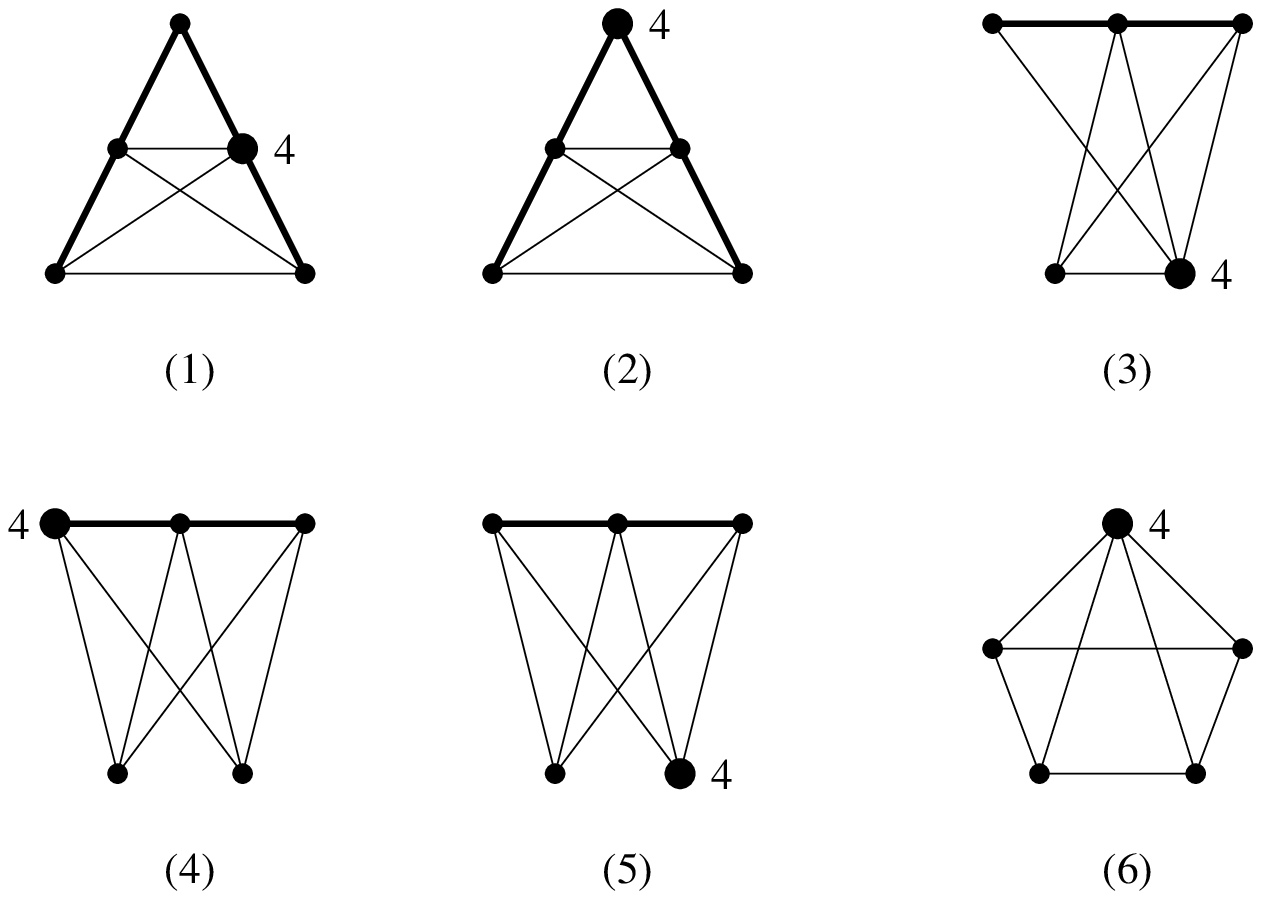} & 63344 & 18 \\ 
\hline
$[2^{10} 3^6]$   & \tabfig{3cm}{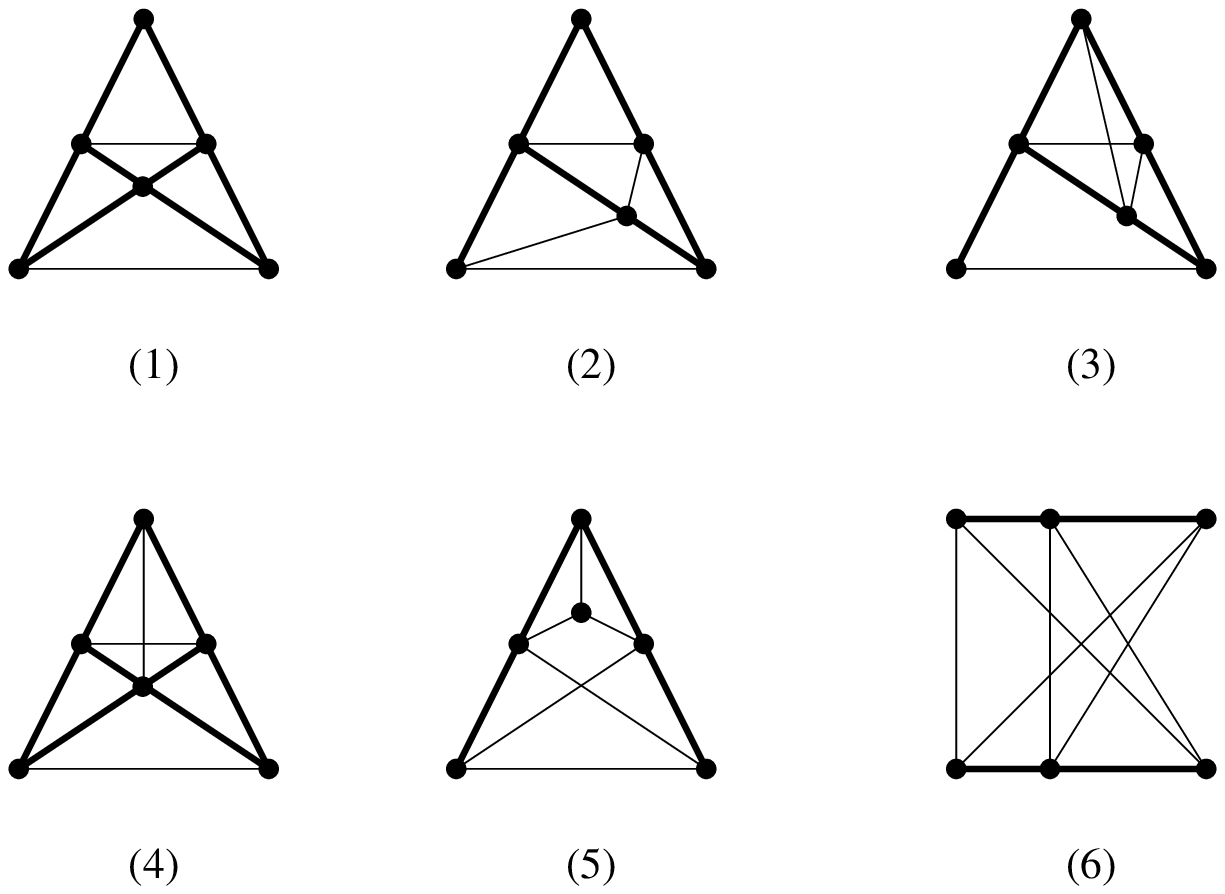} & 88608 & 19 \\ 
\hline
\end{tabular}
\newpage
\begin{tabular}{|c|c|c|c|}
\hline
Signature & Possible \ lattices & $\equiv$ classes & Similarity\ classes\\
\hline\hline
$[2^{13} 3^1 4^2]$ & \tabfig{1cm}{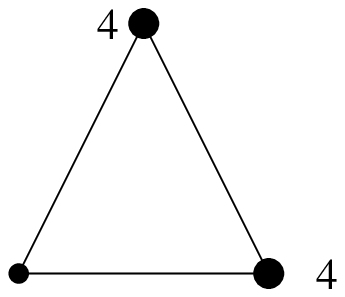} & 21984 & 4\\
\hline
$[2^{13} 3^3 4^1]$ & \tabfig{3cm}{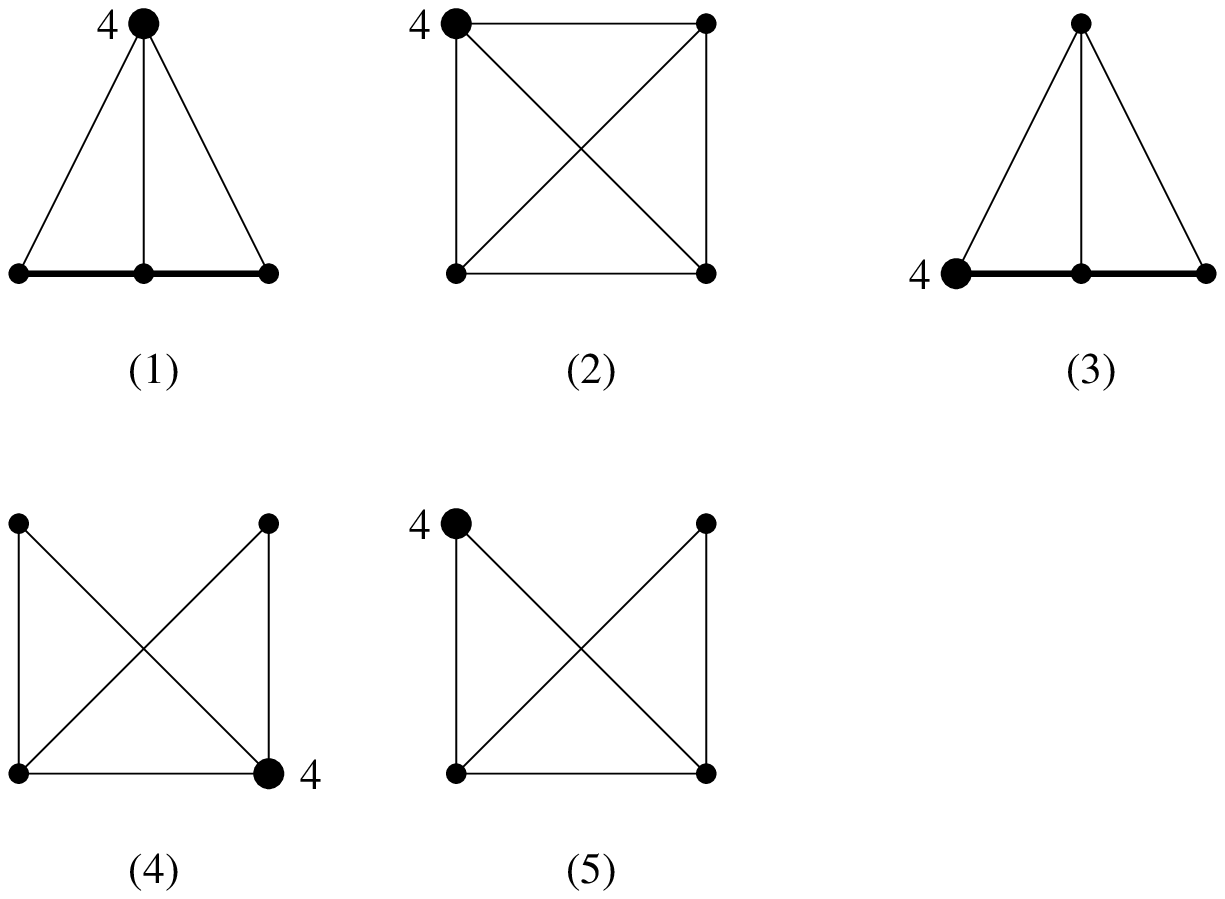}  & 354880 & 22 \\ 
\hline
$[2^{13} 3^5]$   & \tabfig{4.5cm}{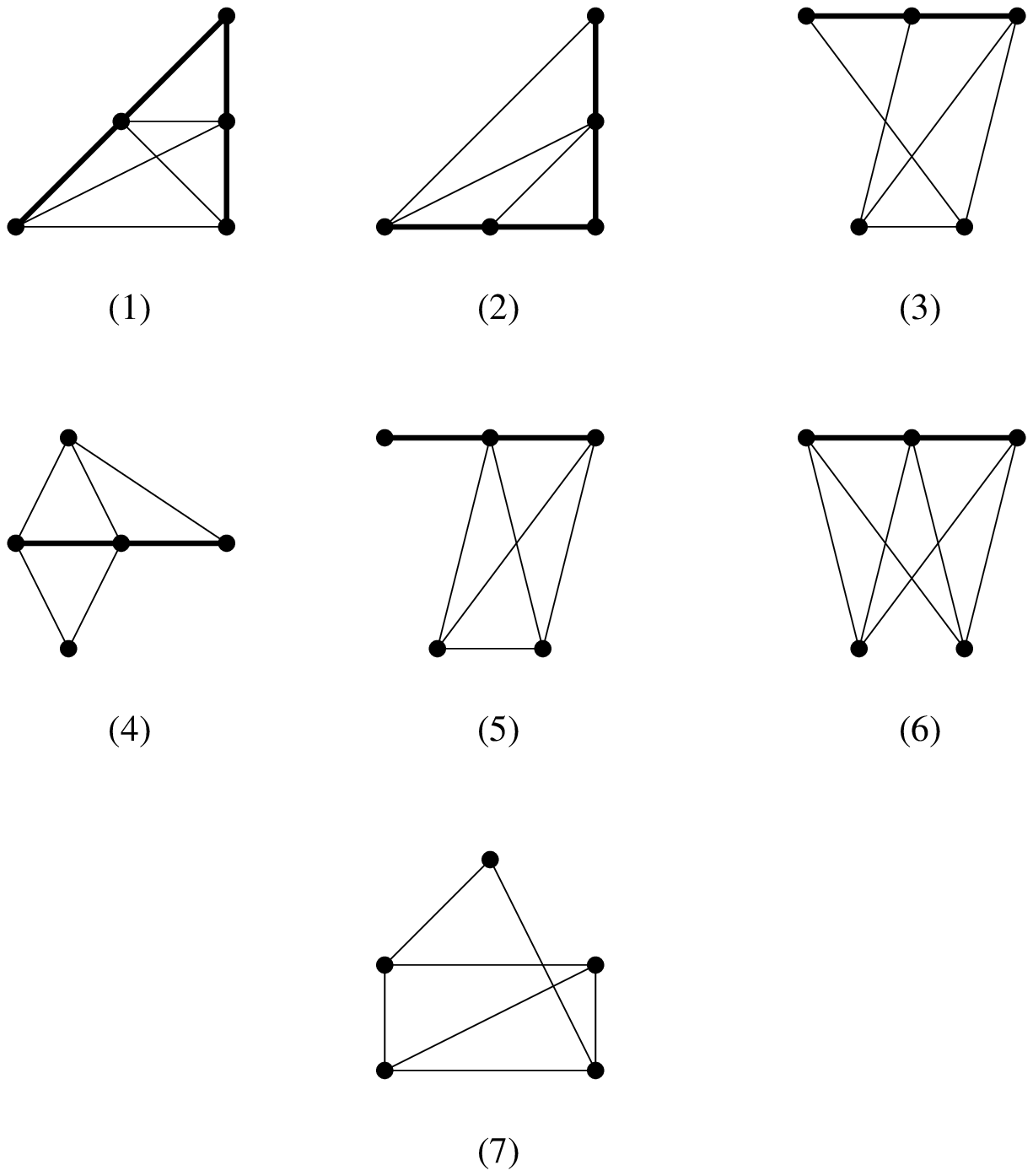} & 776064 & 31 \\ 
\hline
$[2^{16} 3^2 4^1]$   & \tabfig{3cm}{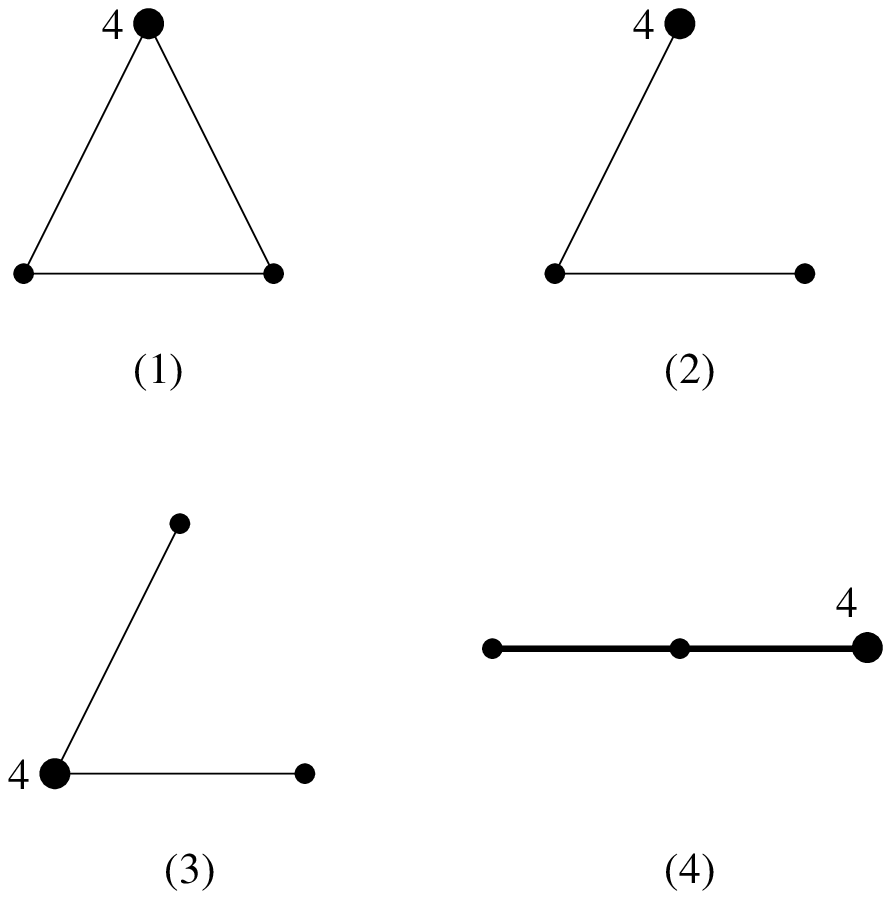} & 884864 & 7 \\
\hline
$[2^{16} 3^4]$   & \tabfig{3cm}{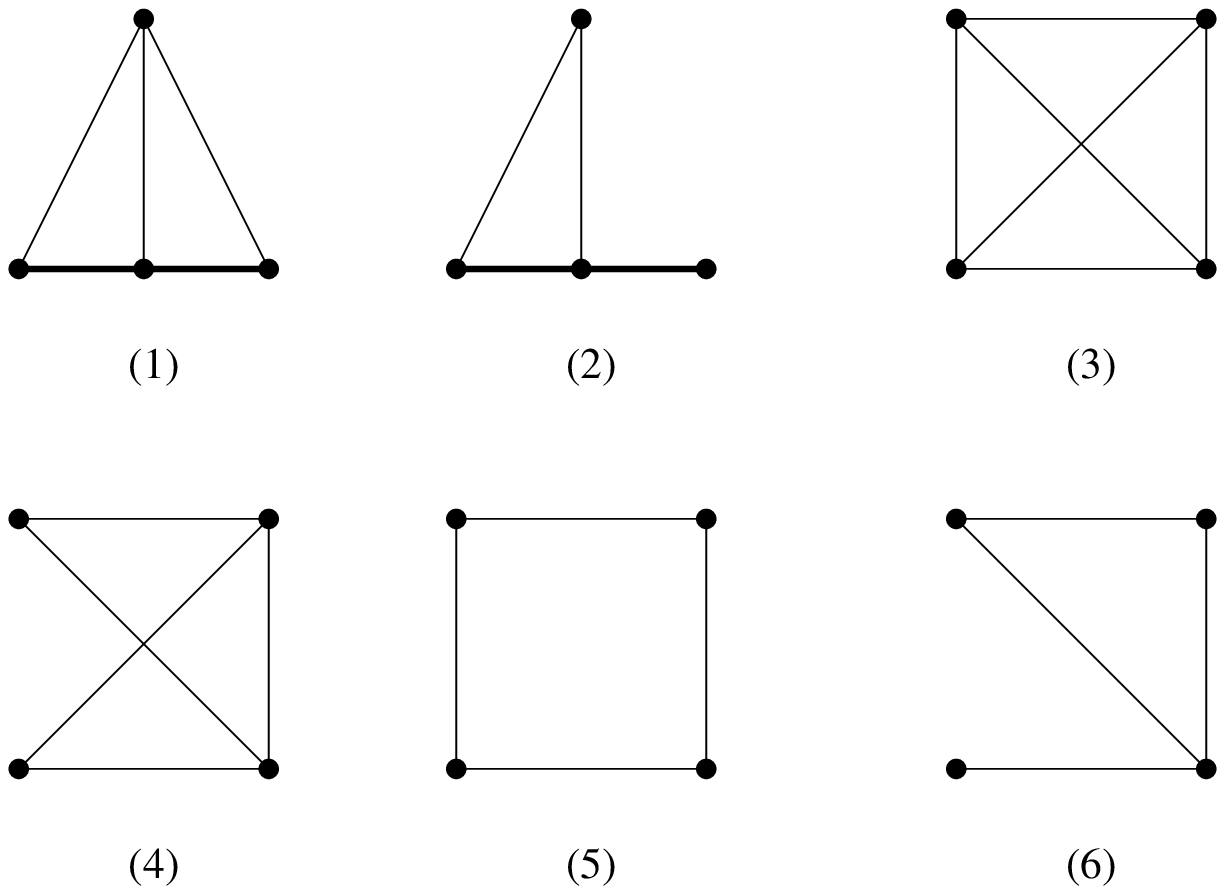} & 3317776 & 18 \\
\hline
$[2^{19} 3^3]$   & \tabfig{3cm}{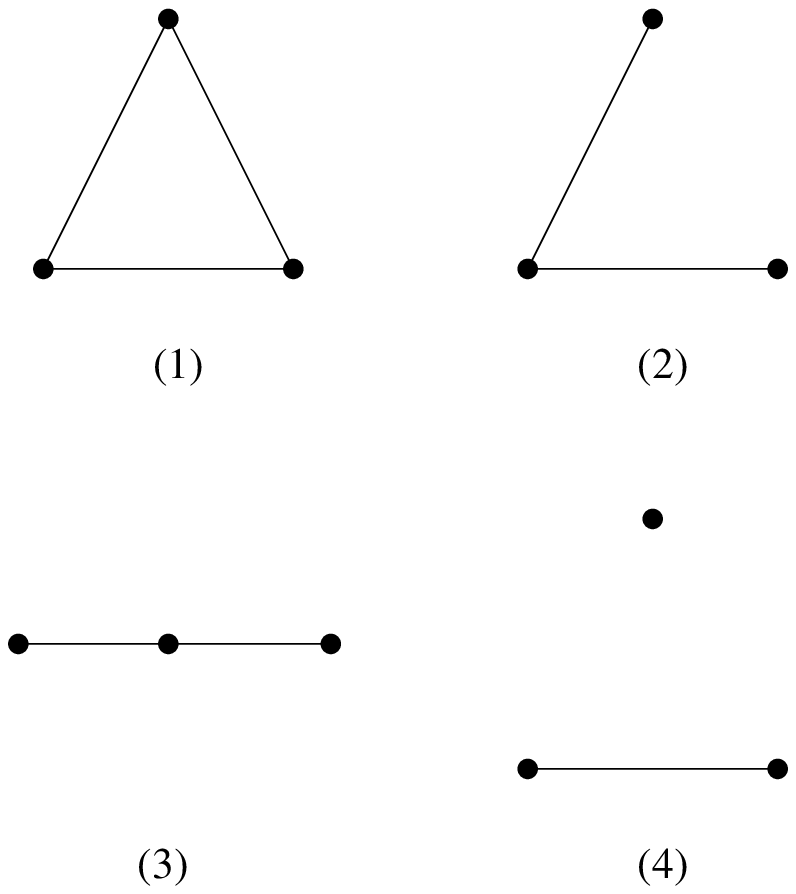} & 7451568 & 5 \\ 
\hline
\end{tabular}
\end{center}

For each signature (here, and in the previous tables), the algorithm described
in Subsection \ref{enumss} produces a list of representatives for the similarity classes. After comparing their incidence lattices (Subsection \ref{latss}), we compute their fundamental groups and compare each pair with the same lattice.
In all the cases we found that the affine and projective fundamental groups
are isomorphic, hence we have the following:
\begin{thm}
For real arrangements of up to $8$ lines, the incidence lattice determines
the affine and projective fundamental groups.
\end{thm}

In two cases we did not have to compare all the groups directly.
For the signature $S = [2^{16} 3^2 4^1]$, in cases (2)--(4) the
graph of multiple points has no cycles, so the projective group
depends only on the lattice by Theorem \ref{Fan}. For the affine
group, Theorem \ref{Fan-} settles case (4), and in cases (2)--(3)
we did compare the appropriate groups. For case (1) there is a
simple line, and then Theorem \ref{simpline} treats the projective
groups, and the special case of Proposition \ref{OS} the affine
groups.

The other signature for which the groups were not compared is $S
= [2^{19}3^3]$ (here we found the largest space $\Wmodequiv$, of
size $7451568$). For lattices (1)--(3), Theorem \ref{simpline}
reduces to seven lines for the projective groups and the special
case of Proposition \ref{OS} for the  affine groups. For case (4)
we use Theorems \ref{Fan} and \ref{Fan-} (alternatively, one can
compute the groups using Propositions \ref{OS} and \ref{Fan_2la}
since the graph of multiple points is disconnected).
As a matter of fact, 
in lattices (2)--(4) there is a unique similarity class, so there 
was nothing to check.

The bottle neck of the algorithm described in Subsection
\ref{enumss} is not so much the running time, but the memory
required to collect the representatives of $\Wmodequiv$. The
extreme cases of $\ell = 8$ lines needed enough memory to upset
our system manager, so we did not continue to $\ell = 9$.

\section*{Acknowledgments}
We thank Sarah Rees 
for her kind help with the \software{testisom} package. Some of
the most stubborn group-comparison problems were solved by her
laptop.


\begin{\bib}{10}
\bibitem[Ar]{Ar} Arvola,~W., {\it The fundamental group of the complement of an
    arrangement of complex hyperplanes}, Topology {\bf 31}, 757--766 (1992).
\bibitem[CS]{CS} Cohen,~D.~C. and Suciu,~A.~I., {\it The braid monodromy of
   plane algebraic curves and hyperplane arrangements}, Comment. Math.
   Helvetici {\bf 72}, 285--315 (1997).
\bibitem[Co]{Co} Cordovil,~R., {\it The fundamental group of the complement of the complexification
   of a real arrangement of hyperplanes}, Adv. Appl. Math. {\bf 21}, 481--498 (1998).
\bibitem[Cox]{Cox} Coxeter,~H.~S.~M., {\it Self Dual Configurations and Regular Graphs}, Bull. A.M.S {\bf 56}, 413--455 (1950).
\bibitem[EHR]{EHR} Epstein,~D.~B.~A., Holt,~D.~H. and Rees,~S.~E., {\it The use of
  Knuth-Bendix methods to solve the word problem in automatic groups},
  J. Symbolic Computation {\bf 12}, 397--414 (1991).
\bibitem[Fa1]{Fa1} Fan,~K.~M., {\it Position of singularities and fundamental group of
the complement of a union of lines}, Proc. Amer. Math. Soc. {\bf 124}(11),
3299--3303 (1996).
\bibitem[Fa2]{Fa2} Fan,~K.~M., {\it Direct product of free groups as the fundamental group
 of the complement of a union of lines}, Michigan Math. J. {\bf 44}(2) 283--291, (1997).
\bibitem[GaTe]{GaTe} Garber,~D. and Teicher,~M., {\it The fundamental group's
  structure of the complement of some configurations of real line
  arrangements}, Complex Analysis and Algebraic Geometry,
  edited by T.~Peternell and F.-O.~Schreyer, de Gruyter, 173--223 (2000).
\bibitem[GKT]{GKT} Garber,~D., Kaplan,~S. and Teicher,~M., {\it A New Algorithm for
   Solving the Word Problem in Braid Groups}, Adv. Math., to appear.
\bibitem[GTV]{GTV} Garber,~D., Teicher,~M. and Vishne,~U., {\it Classes of wiring diagrams
   and their invariants}, submitted.
\bibitem[Ga]{Gthesis} Garber,~D., Ph.D. dissertation, Bar-Ilan University, Israel, 2001.
\bibitem[GJ]{GJ} Garey,~M.~R. and Johnson,~D.~S., {\it Computers and intractability.
   A guide to the theory of NP-completeness}, W.H. Freeman (1979).
\bibitem[Go]{Go} Goodman,~J.~E., {\it Proof of a conjecture of Burr, Gr\"unbaum and
Sloane}, Discrete Math. {\bf 32}, 27--35 (1980).
\bibitem[GP]{GP} Goodman,~J.~E. and Pollack,~R., {\it Allowable sequences
   and ordered types in discrete and computational geometry}, in: New trends in
   discrete and computational geometry, edited by J. Pach, Springer-Verlag,
   103--134 (1993).
\bibitem[HR]{HR} Holt,~D.~F. and Rees,~S.~E., {\it The isomorphism problem for
   finitely presented groups in Groups}, Combinatorics and Geometry,
   London Math. Soc. Lecture Note Ser. {\bf 165},  459--475 (1992).
\bibitem[MH]{MH} Merkle,~R.~C. and Hellman,~M.~E., {\it On the security of multiple
   encryption}, Comm. ACM {\bf 24:7}, 465--467 (1981).
\bibitem[MoTe]{MoTe1} Moishezon,~B.~G. and  Teicher,~M., {\it Braid Group
   Technique in Complex Geometry I, Line Arrangements in $\C\PP ^2$},
   Contemporary Math. {\bf 78}, 425--555 (1988).
\bibitem[OkSa]{OkSa} Oka,~M. and Sakamoto,~K., {\it Product theorem of the
   fundamental group of a reducible curve}, J. Math. Soc. Japan, {\bf 30}(4),
   599--602 (1978).
\bibitem[OrSo]{OrSo} Orlik,~P. and Solomon,~L., {\it Combinatorics and Topology of Complements of Hyperplanes}, Invent. Math. {\bf 56}(2), 167--189 (1980).
\bibitem[OrT]{OT} Orlik,~P. and Terao,~H., {\it Arrangements of Hyperplanes},
   Grundlehren {\bf 300}, Springer-Verlag (1992).
\bibitem[Ra]{Ra} Randell,~R., {\it The fundamental group of the complement of a
   union of complex hyperplanes}, Invent. Math. {\bf 69}, 103--108 (1982).
   {\it Correction}, Invent. Math. {\bf 80}, 467--468 (1985).
\bibitem[Ro]{Rotman} Rotman,~J.~J.,
``An Introduction to the Theory of Groups'',
Springer-Verlag, 1991.

\bibitem[Ry]{Ry} Rybnikov,~G., {\it On the fundamental group of the
   complement of a complex hyperplane arrangement}, preprint (1994)
   [math.AG/9805056].
\bibitem[Sa]{Sa} Salvetti,~M., {\it Topology of the complement of real hyperplanes in
   $\C^N$}, Invent. Math. {\bf 88}, 603--618 (1987).
\bibitem[vK]{VK} van~Kampen,~E.~R., {\it On the fundamental group of
   an algebraic curve}, Amer. J. Math. {\bf 55}, 255--260 (1933).
\bibitem[W]{White} White,~H.~S., {\it The plane figures of seven real lines}, Bull. A.M.S. {\bf 38} 59--65 (1932).
\end{\bib}

\end{document}